\documentclass[11pt]{article}
\usepackage[utf8]{inputenc}
\usepackage[english]{babel}
\usepackage{amsthm}
\usepackage{amsmath}
\usepackage{graphicx}
\usepackage{tabularx}
\usepackage{hyperref}
\usepackage[linesnumbered,ruled,vlined]{algorithm2e}
\usepackage[toc,page]{appendix}
\usepackage{geometry}
\usepackage{bm}
\usepackage{mathtools}

\usepackage{float}
\usepackage{subcaption}

\usepackage[justification=centering, labelfont=bf ]{caption}
\newcommand{\captions}[1]{\linespread{1}\caption{\small{#1}}}
\newcommand{\floor}[1]{\left\lfloor #1 \right\rfloor}

\SetCommentSty{mycommfont}

\SetKwInput{KwInput}{Input}                % Set the Input
\SetKwInput{KwOutput}{Output}              % set the Output

\usepackage{custom_formats}
\usepackage{custom_ian}
\theoremstyle{definition}

\usepackage{setspace}
\providecommand{\keywords}[1]
{
  \small	
  \textit{Key words:} #1
}

\usepackage{titlesec}
\titleformat{\subsubsection}[runin]% runin puts it in the same paragraph
        {\normalfont\bfseries}% formatting commands to apply to the whole heading
        {\thesubsubsection}% the label and number
        {0.5em}% space between label/number and subsection title
        {}% formatting commands applied just to subsection title
\DeclareMathOperator*{\argmax}{argmax} % thin space, limits underneath in displays
\usepackage{authblk}
\title{A Markov process approach to untangling intention versus execution in tennis}
\author[1]{Timothy C.Y. Chan}
\author[2]{Douglas S. Fearing}
\author[1]{Craig Fernandes\footnote{Corresponding author; Email address: craig.fernandes@mail.utoronto.ca (Craig Fernandes)}}
\author[2]{Stephanie Kovalchik}
\affil[1]{Department of Mechanical and Industrial Engineering, University of Toronto, Toronto, Ontario, Canada}
\affil[2]{Zelus Analytics, Austin, Texas, USA}
\date{}             

\begin{document}

\vspace{-4cm}
\maketitle

% \hline

\vspace{-1.5cm}

\begin{abstract}
Value functions are used in sports applications to determine the optimal action players should employ. However, most literature implicitly assumes that the player can perform the prescribed action with known and fixed probability of success. The effect of varying this probability or, equivalently, ``execution error'' in implementing an action (e.g., hitting a tennis ball to a specific location on the court) on the design of optimal strategies, has received limited attention. In this paper, we develop a novel modeling framework based on Markov reward processes and Markov decision processes to investigate how execution error impacts a player's value function and strategy in tennis. We power our models with hundreds of millions of simulated tennis shots with 3D ball and 2D player tracking data. We find that optimal shot selection strategies in tennis become more conservative as execution error grows, and that having perfect execution with the empirical shot selection strategy is roughly equivalent to choosing one or two optimal shots with average execution error. We find that execution error on backhand shots is more costly than on forehand shots, and that optimal shot selection on a serve return is more valuable than on any other shot, over all values of execution error.

%We determine the probability of winning the point starting from a given state, given various shot selection strategies and levels of execution error. 

%Solving our models shows that erring on backhand shots is more costly than on forehand shots. We also show that players should only be aggressive in their shot selection if they have a very low level of execution error. Lastly, we motivate the importance of practicing and exerting the most energy on serve returns, as this has the highest incremental gain in value if played optimally.}

\end{abstract}


\vspace{1cm}
% \hline

\newpage

\section{Introduction}
\label{sec:Intro}

\noindent 
A major stream in sports analytics research seeks to attribute a ``value” to game states in order to assess performance or determine optimal tactical strategies. Most of these models implicitly assume that the athlete or player can perform the prescribed optimal action with a known and fixed probability of success. As an extreme example, consider chess. Given any strategy, a chess player will be able to execute the moves suggested by this strategy in any game state with certainty. There should be no error in the execution of that strategy if it is fully known to the player. The challenge with prescribing optimal actions in other sports, however, is that even when an athlete knows the ``optimal” action, he or she may not be able to execute the action perfectly. For example, an optimal shot in tennis may be to aim at the baseline corner or hit a drop shot, but the execution of these shots is subject to error. Therefore, an optimal action assuming the player can execute perfectly may be different than an optimal action accounting for uncertainty in execution. Moreover, different magnitudes of uncertainty in execution may lead to different optimal actions. In the remainder of this paper, we will refer to uncertainty in the execution of an action as \emph{execution error}.

Previous research has typically accounted for execution error implicitly. For example, uncertainty in the execution of a tennis shot can be captured indirectly via the state transition probabilities. Consider the state being described by the players' locations, the type of shot, and the current score. Now suppose the player chooses an aggressive shot aimed at the sideline. Compared to a conservative shot aimed down the middle, the aggressive shot would likely be associated with greater execution error, which can be captured by an increased probability of the state transitioning to losing the point (e.g., hitting the ball out of bounds). % Now suppose the player who is about to hit the shot chooses a more aggressive shot aimed at the deuce sideline, in response to an incoming shot that is arriving with pace and at a wide angle. Given the difficulty of the incoming shot and the decision to be aggressive, suppose the aggressive shot is prone to higher execution error. Then compared to the same aggressive shot when the incoming shot is less difficult or compared to a conservative shot aimed at the center of the court, the higher execution error associated with a difficult aggressive shot can be captured via a higher transition probability to a state where the player immediately loses the point. %
However, by modeling the execution error of a shot directly (e.g., aiming for a specific location in the court, and the ball landing one foot to the right of that location), we can build a model with more granular control over the probability distribution of state transitions, in order to untangle the impact of intention versus execution in value generation. 

% \sk{This extended example seems overly complex. Could we simply say, "Suppose aggressive shots are associated with greater error. This could be reflected in transition states by an increased probability of transitioning from an aggressive shot to the loss of a point." Or something along these lines?} \cf{will cover}

% probs in state transition
% probs in action execution -> prob in state transition

% 1. we can control error at a much more granular level 
% 2. investigate how this error affects value / strat etc

%To address this challenge, we seek to build optimal policies for players that consider their inherent inaccuracy, or, what we denote, their \emph{execution error}. In this approach, we may recommend a ``strictly sub-optimal shot” that would actually achieve a better result when a particular player's execution error is accounted for.This will lead to more accurate representations of modeling tennis, and allow us to devise the right strategy. This approach also allows us to have different levels of uncertainty in realizing a particular action, dependent on the player we are trying to model. 

In this paper, we develop a modeling framework to investigate how execution error impacts strategy. With direct control over execution error, we can investigate how particular magnitudes and types of execution error affect the value of different game states and the selection of optimal actions. We illustrate our framework within the context of tennis, which conveniently provides a concrete, two-dimensional action space to visualize errors in executing particular shots. In our modeling framework, optimality will be measured with respect to a player winning the point, which is compatible with a rational player's real-world objective. Our modeling framework uses a dynamic programming-based value function. We develop two models, an infinite horizon Markov reward process (MRP) and an infinite horizon Markov decision process (MDP) \citep{Puterman1994}. The MRP is used to compute the value of particular shot selection strategies, while the MDP is used to find an optimal strategy, each for a particular level of execution error.

Our high level approach is as follows. We begin by defining the state space to be the starting locations of both players when the ball is struck and the type of shot (serve, serve return, or rally shot). The action space is a partition of the court into locations where a player would aim their shot -- that is, the player's \emph{intention}. To model the realization of a particular intended shot, we fit bivariate Gaussian distributions on the landing locations of hundreds of millions of simulated tennis shots \citep{kovalchik2020}. The covariance matrices of these distributions model the execution error associated with each intended shot location. We simulate points using these distributions to build state transition probabilities, which are used in our MRP and MDP models. By varying the covariance matrices, we can adjust the modeled execution error at a very fine level, and use the MRP and MDP models to study the efficacy of particular shot selection strategies or to determine an optimal strategy. %We define optimal policy as a mapping between what is the best shot to take in every state that maximizes a player's chance of winning the point. 

This paper makes the following contributions:
\begin{enumerate}
    \item We provide, for the first time in the literature, a general modeling framework to study execution error in professional sports. The advantage of this framework is that we have full control over the level and type of execution error modeled by manipulating the distributions defining the uncertainty in action outcomes. This control facilitates the quantification of the impact of execution error on value and optimal strategies.
    \item In the context of tennis, we demonstrate several numerical insights on the relationship between execution error and value, namely:
    \begin{enumerate}
        \item Execution error on backhand rally shots are more detrimental to value compared to errors on forehand rally shots.
        \item Execution error on serves do not greatly impact value.
        \item As execution error increases, the optimal strategies for players become more ``conservative", with more shots aimed toward the middle of the court.
        \item If players could only choose one shot to exert more energy and strive for optimality, it is most valuable to do so on the serve return.
    \end{enumerate}
\end{enumerate}

\section{Related Literature}
\label{sec:LitReview}

\noindent 
This paper relates to three bodies of sports research: (1) dynamic programming and Markov processes, (2) shot and point predictions in tennis and (3) optimizing strategies for execution error.
%\tcyc{i think we should stick with american english throughout: modelling -> modeling. or we can have a discussion first about journal target before we make changes (i have already edited intro to make it single l, but will not make changes from here on out until we decide.)} \cf{Changed the rest to single l.}

Dynamic programming has been extensively explored within sports to solve tactical decision problems and compare teams and players. \cite{Fry2007} defined a deterministic dynamic program to model the series of decisions coaches must make during a draft to maximize their total draft value.  \cite{Broadie2008,Broadie2012} developed a novel metric, \textit{strokes gained}, that identifies the strengths and weaknesses of a golfer's game by building a value function in golf through dynamic programming. \cite{Chan2021} developed a similar metric, \emph{points gained}, in American football and provided a general theoretical justification for the use of dynamic programming-base value functions in sports. \cite{Hirotsu2019} modeled baseball as an MDP to quantify the value of a sacrifice bunt, discovering it to be more beneficial than previously thought. Within tennis, \citet{Clarke2015} used dynamic programming to determine the optimal policy of when a player should challenge an umpire's call and \cite{Chan2016, Chan2018} used Markov decision processes and Bayesian regression to develop a tennis handicap system.% \tcyc{add raghav's second paper on this topic within the same cite} \cf{- done}. 

In tennis, there has been a significant focus on developing predictions on individual shots and points. \citet{Wei2013} used Bayesian networks to compute an individual shot outcome prediction. \citet{Wei2016a} and \citet{Wei2016b} use statistical models to discover patterns in shot sequences to predict the next likely shot and compute the probability of winning the point. Inspired by \citet{Cervone2016}, \citet{Floyd2020} utilizes multi-resolution stochastic modeling to develop \textit{expected shot value}, a metric that assigns shot value to individual shots throughout a point. \citet{kovalchik2020} also estimates expected shot value through a mixture modeling approach of both ball and player trajectories. Lastly, \citet{Terroba2013} developed a detailed MDP with states that include both players' starting locations and the incoming shot difficulty, and an action set that specifies the stroke, shot direction and shot type. They develop video recognition algorithms to build a tennis data set and used Monte Carlo tree search methods to solve the MDP and compute optimal shot selection strategies. 

To date, execution error has been relatively unexplored and none of the papers discussed so far have considered it in their analyses. The few studies that have examined execution error include \citet{Ahmad2016} and \citet{Yee2016}, who developed models to devise optimal shot strategy in curling. They incorporate an execution error noise term by using fat-tailed t-distributions around the intended shot to calibrate the shot outcomes on par with Olympic level players. \citet{Tibshirani2011} fit bivariate Gaussian distributions on dart throws to model the execution error around a specified target. \citet{Haugh2020} built upon this work by feeding these distributions into a dynamic program and dynamic zero-sum game to develop an optimal strategy of where to aim their dart, depending on a player's specific execution error, the current score, and the opponent's strategy. 

Our paper adds to the growing literature on execution error by introducing it to the tennis context and, more importantly, developing a framework that treats execution error as a ``controllable lever" that can be used to explore how execution error influences value creation. %This is done in a robust and replicable \tcyc{do robust and replicable have a precise meaning here? are the prev approaches not robust/replicable?} \cf{robust/replicable are definitely misnomers. I am trying to highlight the ``lever" on epsilon that we have, as past literature didn't play with magnitude/type of execution error, just added an error term in the execution of actions.} within the MRP and MDP transition matrix framework.

\section{Data}
\label{sec:Data}

\noindent

Our data consists of simulated tennis shots generated through the VON CRAMM framework, created by \citet{kovalchik2020}. Their framework provides us with accurate 3D ball and 2D player trajectories that are in strong agreement with real grand slam tennis data. This is done through an infinite Bayesian Gaussian mixture model that provides a generative distribution for ball and player trajectories. This model is fit using highly scalable variational inference methods trained on 125,000 men's shots and 80,000 women's shots sampled from past Australian Opens. Additionally, the framework allows for conditional distributions for the generative model to be built that can project all possible future paths for the ball and players. Lastly, the shot outcome of ``win'', ``error'' or ``in-play'' is predicted from hierarchical generalized additive models with non-linear spatial effects.

\begin{figure}[htbp]
    \centering
    \includegraphics[width=12.5cm]{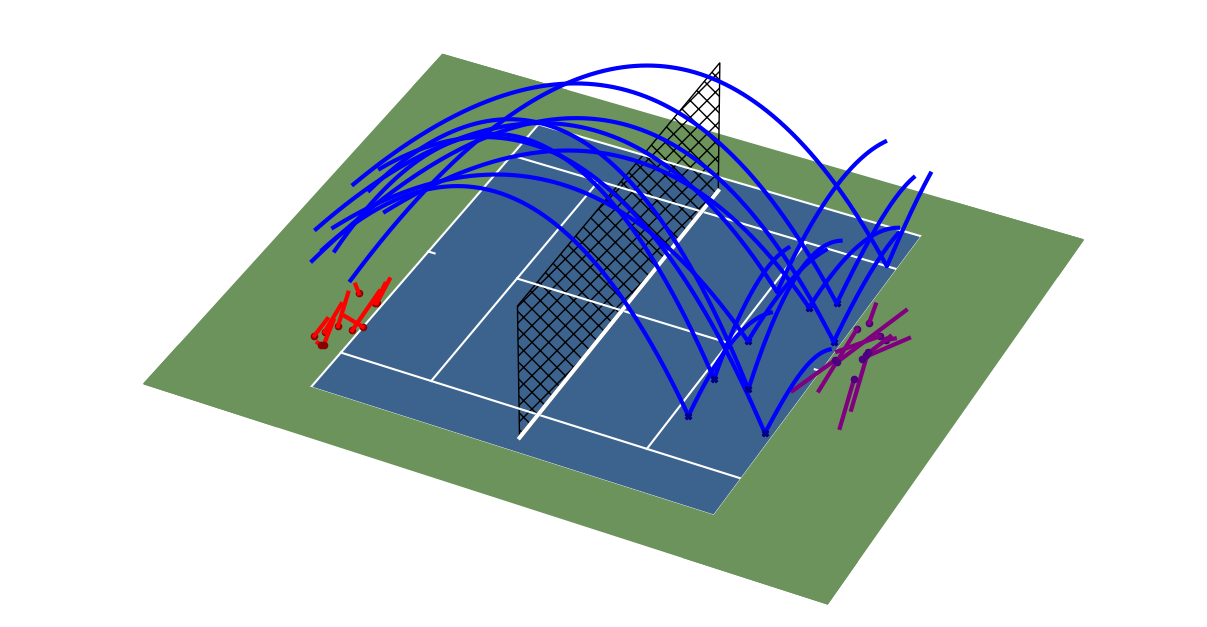}
    \caption{Visualization of 10 rally shots generated from the VON CRAMM framework. The shots were conditioned such that at the moment of impact, the player who hit the ball (red) was standing near the deuce corner and the receiving player was standing at the center mark. The red and purple dots indicate the starting positions of the two players.}
    \label{fig:data-vis}
\end{figure}

The VON CRAMM framework offers flexibility in generating tennis shots and modeling the game. First, we can generate millions of realistic tennis shots, each associated with a prediction of whether it won the point or not. Simulated shots will be drawn from a distribution $G$, which we refer to as the shot generator. Second, and more importantly, we can condition on features such as the locations of the players and the landing location of the shot when generating the trajectory of a shot. When we use $G$ to generate a shot, $w$, conditional of features $\xi $, we denote this by $w \sim G \: | \: \xi$. %Features that we can condition on include the locations of the players and the landing location of the shot. 
%For example, if we wanted to model a shot where the player who is about to hit the shot is standing near the net and the receiver is standing at the baseline, and the ball is aimed toward the sideline, that is easily achievable assuming we have a distribution that relates the actual ball landing location to the location where the ball is aimed. 
Figure \ref{fig:data-vis} provides a visualization of several rally shots generated by VON CRAMM. %In the remainder of this paper, we will refer to the player who is about to hit the shot as the \emph{impact player}. 
% \tcyc{[commented out impact player since we only used it 2-3 times]}

A limitation of this framework is that each shot is generated independently, so we cannot directly model a particular sequence of shots. However, we can achieve this indirectly by making use of the conditioning ability of VON CRAMM. For example, we can generate a serve, examine the terminal features of that shot (i.e., where the players are standing when the receiver is about to strike the ball), then generate a likely return shot conditioned on those features, and so on. In this fashion, we can ``stitch" together a likely sequence of shots to model a point. These player transitions between shots will serve as the underlying transition probabilities we use in our MRP and MDP models. %\tcyc{do you mean prob transitions specifically?} \cf{yeah - edited accordingly}

Finally, we note that the VON CRAMM framework is actively used by Tennis Australia. It enabled the creation of a novel tennis metric, \textit{serve value}, which measures the quality of a player's serve at the Australian Open in Melbourne \citep{cohen2020}.

\section{Markov Process Models}
\label{sec:Method}

\noindent
To develop value functions and optimal policies in tennis, we first need to model a point in tennis, which we do within the framework of a Markov decision process. In this section, we describe the components of our model, demonstrate how to incorporate execution error as an adjustable parameter, and formalize policy evaluation and policy optimization considering this error, through a Markov reward process model and a Markov decision process model, respectively.

%\tcyc{will return to this intro paragraph after finishing rest of section...}
%To do this, we need to model a player's action set, and how players traverse the court for any particular chosen action. We then seek to artificially introduce execution error and model how this affects the player transitions, overall shot value and optimal policy. We accomplish this through two models, an MRP and MDP \citep{Puterman1994}, which we now describe in detail. 

\subsection{Model components}

%\tcyc{should we change these to subsubsections?}

\subsubsection{States} \hfill\\ \vspace{-0.5cm}

We define the state as $s = (\sigma_A, \sigma_B, \omega)$, where $\sigma_A$ and $\sigma_B$ refer to the locations of Player A and Player B, respectively, and $\omega$ refers to the type of shot (serve, serve-return or rally). We model everything from the perspective of Player A, meaning that $\sigma_A$, $\sigma_B$ and $\omega$ all correspond to the moment that Player A strikes the ball. By defining the state in this manner, we are not concerned with the  ``state" of the system when Player B strikes the ball, only when the ball returns to Player A's side and is about to be hit again. %Rather, this can be seen as a ``dummy state" from the perspective of Player A.  
Lastly, we define two absorbing states. We denote by $W$ the ``win" state and by $L$ the ``lose" state, corresponding to Player A winning or losing the point, respectively. We define $S$ as the set of all possible states in a point.

\begin{figure}[H]
\centering
\begin{tabular}{cccc}
    \includegraphics[width=14.5cm]{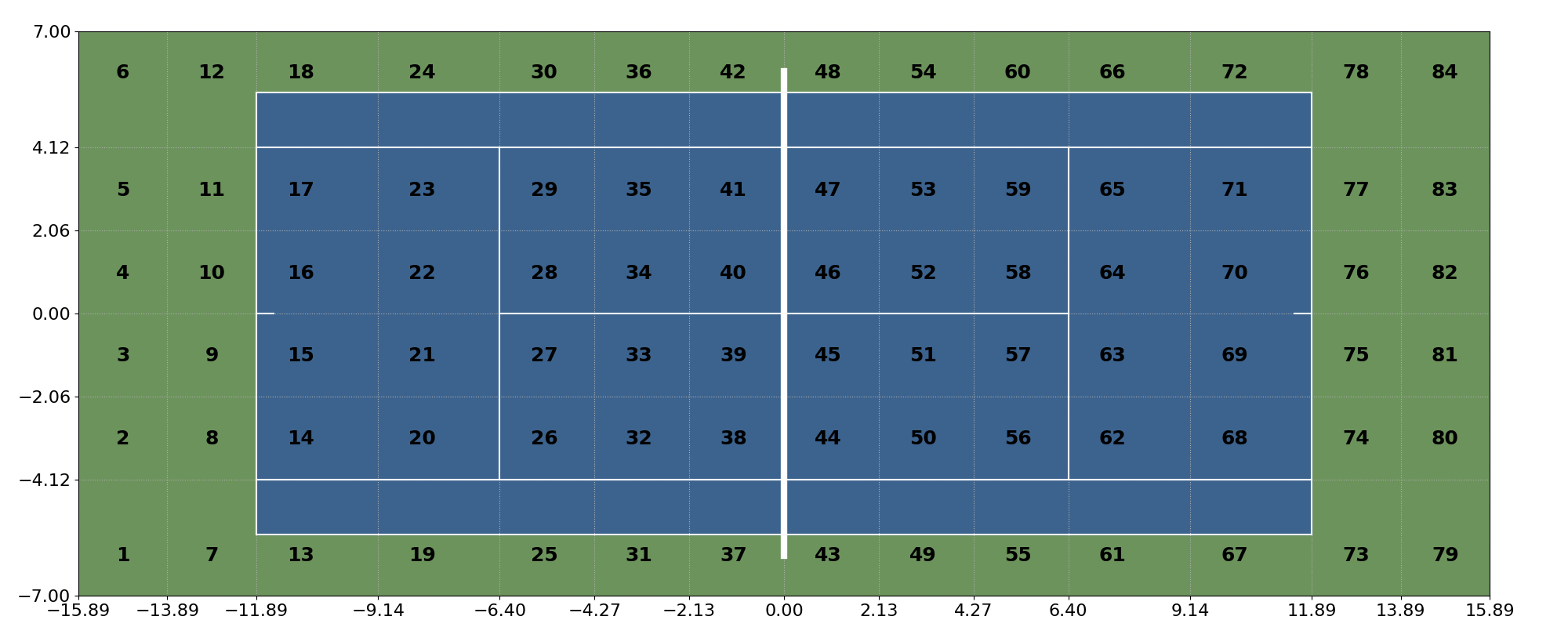} \\
    (a) Discretized locations for the state space, describing the possible values for $\sigma_A$ and $\sigma_B$\\[6pt]
\end{tabular}
\begin{tabular}{cccc}
    \includegraphics[width=14.5cm]{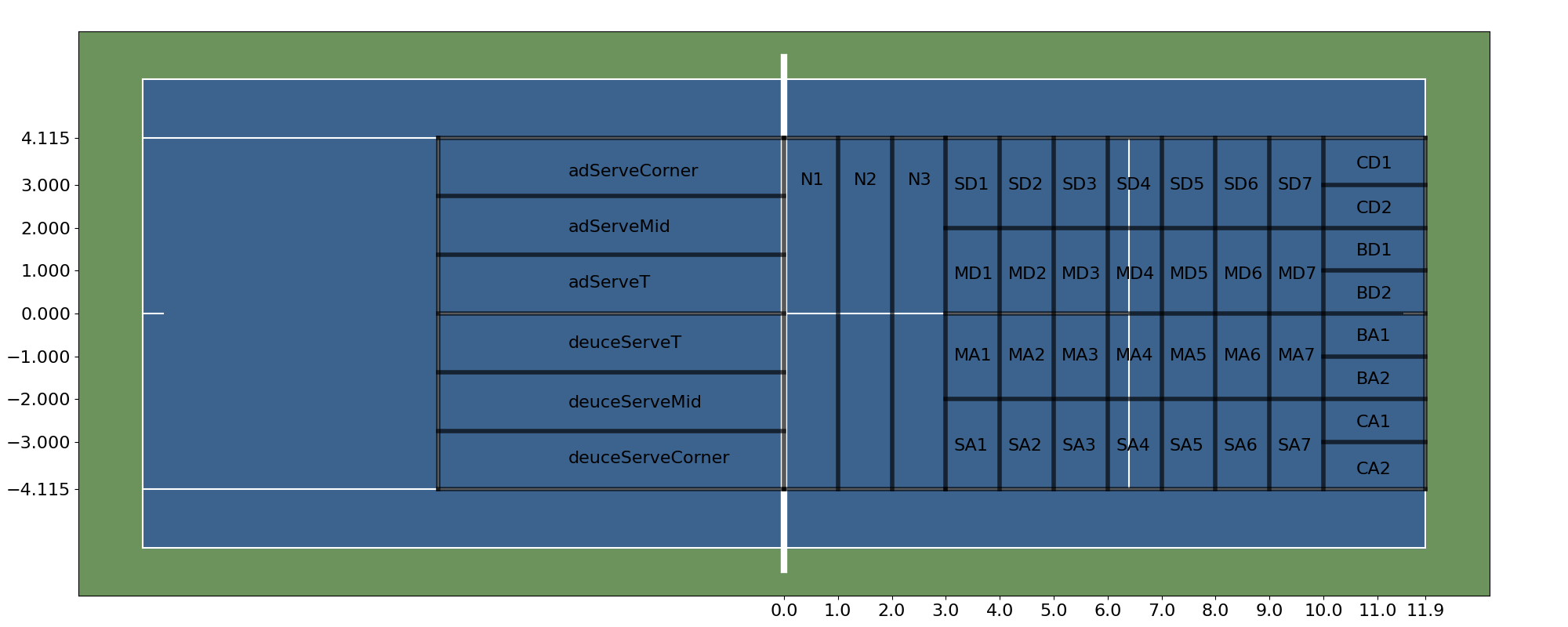} \\
    (b) Action set \\[6pt]
\end{tabular}
\caption{The state space (a) and action set (b) for the MRP and MDP. The values on the x and y axis denote the distance (in meters) from the net and centre line, respectively. The actions on the left side of the court illustrate the actions for serves (serving from the right side of the court). The actions on the right side of the court correspond to serve return and rally shot actions (hitting from the left side of the court). The action abbreviations are Net (N), Sideline Deuce (SD), Middle Deuce (MD), Middle Ad (MA), Sideline Ad (SA), Corner Deuce (CD), Baseline Deuce (BD), Baseline Ad (BA), and Corner Ad (CA).}
\label{fig:states_actions}
\end{figure}

We segment the court into 84 cells of roughly two square meters, 42 on each side of the court (see Figure \ref{fig:states_actions}(a)). The total number of possible transient states is $42 \times 42 \times 3 = 5,292$, since $\sigma_A \in \{1, \ldots, 42\}$, $\sigma_B \in \{43, \ldots, 84\}$, and there are three possible shot types. After removing impossible states (e.g., serving player standing at the net) or very unlikely states (e.g., receiving player standing on the ad side when the other player is about to serve from the deuce side), we are left with 3,600 states for serves, serve returns and rallies (the transient states) and two states for winning or losing the point (the absorbing states). 

%\tcyc{i don't see how you get 72 for serving states. i assume player A stands in only 8 or 9, but Player B could be anywhere on their half of the court. if we start to make realistic judgments about where Player B could stand, then we should also start doing the same for service returns, which would result in fewer than 1764 for service returns...}

%Since serves typically originate from particular areas on the court, we only consider 72 serving states. For serve returns and rallies, we use all combinations for the 42 x 42 grid, resulting in 1764 states for each shot type. In total, we have 3600 states for serves, serve returns and rallies (the transient states) and 2 states for winning or losing the point (the absorbing states). \\

\subsubsection{Actions} \hfill\\ \vspace{-0.5cm}

In each transient state $s$, Player A chooses an action $a$ from the set of permissible actions in that state, denoted $A_s$. The set of all actions is denoted $A \coloneqq \cup_{s \in S} A_s$. Each action refers to a particular region on the court that Player A's shot can be aimed. As such, we refer to actions synonymously as \emph{intentions}, since they refer to the intended location of a shot, which may differ from the realized shot location once execution error is considered. The possible intentions differ for serve versus non-serve shots. As shown in Figure \ref{fig:states_actions}(b), there are six serve actions (three on the deuce side and three on the ad side) and 39 actions for serve return and rally shots. These actions were designed to span the possible locations a ball could be aimed on the court, with larger action regions closer to the net since fewer balls land in those regions. %\cf{Should we mention here how we came up with these 45 action definitions?}

A deterministic \emph{policy} is a function $\pi: S \rightarrow A$ that maps all states to an action. We refer to the corresponding randomized policy as an \emph{intention distribution}. That is, the intention distribution for a given state is a probability distribution over the set of actions. For each intention, there will be an \emph{execution distribution} -- the conditional distribution of landing locations for the ball given the player's chosen intention. The intention and execution distributions are used in the generation of the transition probabilities. % Next, we describe how we fit these distributions using data.

%The probability distribution dictating how often each action is chosen in each state is defined as the \emph{intention distribution}. \tcyc{this is simply a randomized policy, right? should we harmonize the lingo?} It determines which intended shot Player A should choose in any given state. 

%In an ideal world, we would always know where a player \emph{intended} to aim their shot. However, in practice we only know where the ball actually landed. This relationship between intention and execution is what we try to untangle. For every state, a player will have an intention distribution - the distribution indicating what regions of the court they will aim. 

\subsubsection{Transition probabilities} \hfill\\ \vspace{-0.5cm}

We first present the general framework for deriving the transition probabilities, which relies on the intention distribution for each state and the execution distribution for each intention, and then describe how we fit the intention and execution distributions. 

Given the current state $s$ and a chosen intention $a$, we simulate a ball landing location according to that intention's execution distribution. Then we use VON CRAMM to generate Player A's shot conditional on the state $s$ and the landing location. VON CRAMM is used again to generate Player B's return shot with the state being the terminal conditions of Player A's simulated shot. The terminal conditions of Player B's shot correspond to the new state $s'$ where Player A is about to hit the ball again (unless the point has been won or lost within those two shots). For a given state $s$, this process is repeated $N$ times, each time drawing an intention $a$ from the intention distribution. Finally, this process is repeated for all states $s$, each with its own fitted intention distribution. The result of this process is an empirical distribution for $P(s' | s, a)$. 

Note that VON CRAMM doesn't distinguish between generating a first serve and a second serve. Instead, it generates an ``average" serve over the first and second serve characteristics, which is a limitation of our data.  %Although in practice, first and second serves differ in style, we accept this limitation as it doesn't greatly impact our ability to create a modeling framework. 
An implication of this limitation is that if a simulated serve is a fault, we need to determine whether this was a first serve or a second serve to decide whether the serve should be repeated or transition directly to the state $L$. Using a public tennis data set \citep{Sackman2021}, we use Bayes' rule to compute $P(\text{first serve} \: | \: \text{fault}) = P(\text{fault} \: | \: \text{first serve})P(\text{first serve})/ P(\text{fault}) = 0.379\cdot0.725/0.302 = 0.91$. 
%
%\sk{[Is this essentially saying that a fault is 10 times more likely on first serve than second? I have more intuition around $P(fault \: | \: first serve)$ so trying to get my head around what $P(first serve \: | \: fault)$ should be.]} \cf{[Using the data set, I computed $P(first serve \: | \: fault) = \frac{P(fault \: | \: first serve) * P(first serve)}{ P(fault)} = \frac{0.379*0.725}{0.302}$.]}
%
Therefore, when we encounter a fault in a serving state, we encode the probability of a self-transition as $0.91$ and a probability of transitioning to the lose state as $0.09$. Although this simplified approach could lead to more than two serves being modeled, this situation occurs only with probability 0.008, which is negligible for our purposes. 

% \cf{[QUESTION to Stephanie: is there anything else we can add to defend the limitation of first vs second serves not being distinguishable?]}

% \sk{[Without the first and second serve conditioning, what we have is a distribution of all serve types that are 'available' to players and this would tend to be more conservative on first than expected and riskier on second than expected (in the absence of any adjustment). The main difference between 2nd and 1st serves is that players tend to avoid the center of the service box on 1st serve. Would it be possible to define a first and second serve intention region based on this typical pattern on serve?]}

\textbf{Fitting the intention distribution}: %To determine the intentions that make up the intention distribution (the support of the intention distribution), we primarily utilized tennis domain knowledge (see Figure \ref{fig:states_actions}). The only modeling constraint was that the execution distribution for each intention needed to follow a bivariate Gaussian distribution as closely as possible. There is a particular action set for ad serves (3 actions), deuce serves (3 actions), and a common action set for serve returns and rally shots (39 actions). The intention distribution only covers the inbound region of the court, as a player would never purposely hit the ball out-of-bounds. 
%\tcyc{the commented paragraph doesn't describe how we fit the intention distribution. this paragraph should say something like, given a state s, we use von cramm to generate $N$ shots. we then build the intention distribution simply based on the empirical count of the number of shots that landed in each action region.} \cf{True - but then maybe what you commented out goes in 4.1.2, because there should be some (1-2 sentences probably) descibing how we delineated each action (i.e., MA1, MA2, ...)}
To determine the intention distribution for a given state $s$, we generate $N$ shots, $w_1, \ldots, w_N$ using VON CRAMM conditioned on $s$ ($w_1, \ldots, w_N \sim G \: | \: s$) and compute the proportion of shots that landed in each action region, normalized by the total number of shots that landed in bound. We denote this distribution as $f_s$, that is $f_s(a) := P(a \: | \: s)$ is the probability of choosing action $a$ in state $s$. %We refer to it as the \emph{empirical intention distribution} since it was estimated from the ``empirical'' (simulated) data generated by VON CRAMM. \tcyc{i ctrl-f and found that we don't use empirical intention distribution much. maybe we should get rid of `empirical'. do we use ``intention distribution'' anywhere in the paper that doesn't actually refer to the ``empirical'' version of the intention dist? if not, then we might as well simplify and just use ``inte dist'' w/o `empirical' everywhere.}
% \tcyc{i got rid of empirical intention dist. we didn't seem to use it much and i'll keep reading to see if it reads ok.}

While execution error is embedded in the shots generated by VON CRAMM, we expect that these errors cancel out to a certain extent so that the overall distribution, conditioned on inbound shots, approximately represents the intention distribution. For example, given a certain intention, some shots aimed there will spill over into neighboring intentions due to execution error. However, shots aimed at neighboring intentions will also spill over into the intention under consideration. Where this approximation is weakest is likely around the sidelines and baseline, where many shots that landed just out of bounds would have been aimed at the closest inbound intention but not be counted in our conditional distribution. To test the sensitivity of this approach, we generated 100,000 shots and noted that only 12\% were out of bounds. Attributing them to the closest inbound intention did not alter the intention distribution much. Thus, we kept our original approach of only using the inbound shots.

\begin{figure}[H]
\centering
\begin{tabular}{cccc}
    \includegraphics[width=12.5cm]{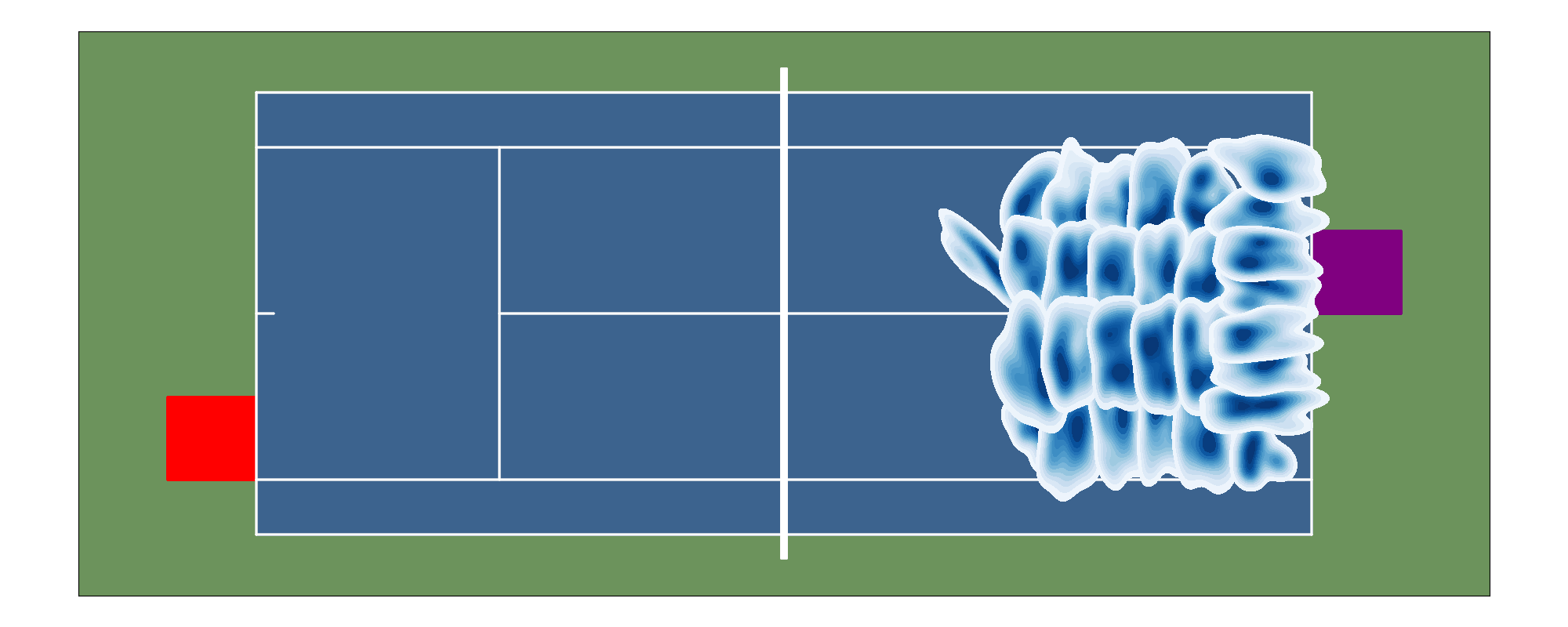} \\
    (a) \\[6pt]
\end{tabular}
\begin{tabular}{cccc}
    \includegraphics[width=12.5cm]{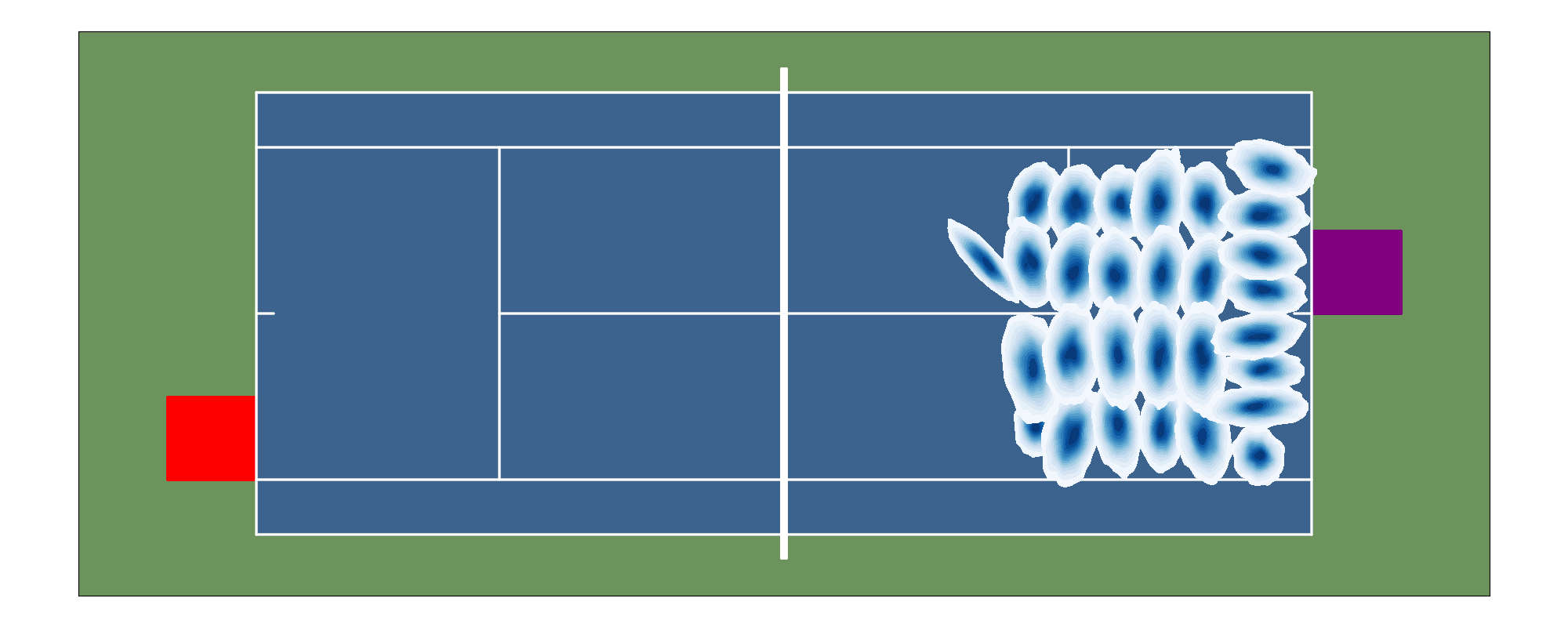} \\
    (b) \\[6pt]
\end{tabular}
\caption{(a) The raw landing location distribution for actions when Player A is in cell 10 (red square), Player B is in cell 81 (purple square) and Player A is hitting a rally shot. (b) The fitted bivariate Gaussians for each action.}
\label{fig:fitted-actions}
\end{figure}

\textbf{Fitting the execution distribution}: We first fit a distribution of ball landing locations for each action that represents a \emph{perfect execution} scenario. Conceptually, we think of perfect execution to mean that the ball's landing location matches the intention. This implies that there may be slight variations in the ball's exact landing location, but not to the point that it will land outside of the particular intention. Given a set of shot landing locations within a particular intention, we fit a bivariate Gaussian distribution and then scale down the covariance matrix so 90\% of the distribution's probability mass lies within the intention. This choice was a balance between having a distribution that represents perfect execution (where 100\% of the shots should land inside the intention), without having too small a covariance matrix. Figure \ref{fig:fitted-actions}(a) shows the raw distribution of landing locations for 1000 simulated shots given the state where Player A is standing in the red cell and Player B is standing in the purple cell. Figure \ref{fig:fitted-actions}(b) shows the fitted bivariate Gaussian distributions representing perfect execution of the landing locations for each intention. For a given state $s$ and given action $a$, the perfect execution distribution associated with this action is written as $\mathcal{N}(\bm{\mu}_{s,a},\bm{\Sigma}_{s,a})$, where $\bm{\mu}_{s,a} \in \Re^2$ and $\bm{\Sigma}_{s,a} \in \Re^{2\times2}$. %\rs{[What is random here? Action a is fixed and the execution is perfect, no? Are you assuming that even perfect execution is not 100\% perfect but only 90\% perfect?]} \cf{[I try to clarify this above]}
%\tcyc{should we write this as dependent on both s and a explicitly? if we only write a and no s, feels like it could mislead that the execution error for each a is the same across all s, which is not the case right?} \cf{yeah I like this idea}

%\tcyc{what i was struggling with the last little while with the writeup was that there were some repeated parts in action as well as transition prob. but i think reorganizing the way i've done above removes redunancy and makes things more clear} \cf{agree - this looks cleaner}

%This corresponds to a player always hitting the ball within their intended action's region. A perfect execution player would never hit the ball out of bounds, so to fit the distributions, we condition only on shots that landed inbound. Then looking at all the shots that landed in a particular region, we compute the mean ($\bm{\mu}$) and covariance ($\bm{\Sigma}$) of those shot's landing locations. This allows us to define the bivariate Gaussian as $\mathcal{N}_2(\bm{\mu}, \bm{\Sigma})$. In order to ensure that samples drawn from this fitted distribution are contained to that particular action's region of the court, we introduce a multiplicative adjustment parameter, $k$, on the covariances to scale down the sizes. The largest adjustment parameter that allows samples to be held within an action's region with at least 90\% consistency was $k=0.7$. This scenario allows us to recapture the intention distribution with the execution distribution.

%We then fit the bivariate Gaussian execution distributions for each intention and scale the covariances down by 0.7. 

%\tcyc{perhaps we should use "empirical intention distribution" as the one derived from data, and "intention distribution" as a generic term for a randomized policy.} 

%To illustrate, let's consider the state where Player A is standing in cell 8, Player B is standing in cell 76 and Player A is about to hit a rally shot, denoted as $s = (8,76,\text{rally})$. 

%Now that we have these two distributions, to generate the transition probabilities, we can follow a five step approach: (1) draw Player A's intention from the intention distribution (2) draw a sample from that intention's corresponding execution distribution to determine a landing location (3) use VON CRAMM to conditionally generate Player A's shot with the parameters from state $s$ and the sampled landing location (4) use VON CRAMM to conditionally generate Player B's return shot based on the ending conditions of Player A's shot (5) record the new state when Player A is about to strike the ball again (unless the point has been won or lost within those two shots). This is repeated for all states.

\textbf{Scaling the execution error:} We now describe how to compute the transition probabilities while modeling different levels of execution error. We define $\epsilon$ to be a controllable parameter that scales the covariance matrices of the bivariate Gaussians. Mathematically, we use $\epsilon$ to model the ``level'' of execution error associated with action $a$ in state $s$ through $\mathcal{N}(\bm{\mu}_{s,a}, \epsilon\bm{\Sigma}_{s,a})$. Note that $\epsilon = 1$ recovers the perfect execution scenario. Larger values of $\epsilon$ lead to more spread in the Gaussian distribution of landing locations, thereby increasing the probability that a shot will land in a location that does not coincide with the intention (i.e., larger execution error). An implicit assumption in this approach is that the intention distribution remains constant for different levels of execution error. This assumption can be relaxed through the following approach for each state: (1) create the intention and execution distributions as before; (2) draw $\floor{N\cdot f_s(a)}$ samples from $\mathcal{N}(\bm{\mu}_{s,a}, \epsilon\bm{\Sigma}_{s,a})$ for each $a$; (3) re-fit the intention distribution based on these new samples. However, these additional steps add computational expense. Thus, for convenience, we keep the same intention distribution for each value of $\epsilon$. We write $P_{\epsilon}(s' \: | \: s,a)$ to denote the transition probability from $s$ to $s'$ when action $a$ is chosen and subject to execution error of $\epsilon$.

%We define $p^\epsilon_{s,s' | a} \coloneqq P(s' \: | \: s,a, \epsilon) $ to denote the probability of transition from state $s$ to $s'$ when action $a$ is used for a player with error level $\epsilon$. 

%wanted to model Player A having ``execution error level $\epsilon$", Player A's shots for a particular intention would follow $\mathcal{N}_2(\bm{\mu}, \bm{\Sigma}\epsilon)$. The greater the error, the larger the covariance. For players with larger execution error levels, the distance between where a player intended to hit the ball and where it actually landed will grow. 

%\tcyc{PAUSED HERE}

The full pseudo-code for generating transition probabilities in the presence of execution error is displayed in Algorithm \ref{algo:TransProbs}. Note that we define $\Bar{s}$ as an auxiliary state that captures the terminal conditions of Player A's shot, which forms the state that Player B's generated return shot is conditioned on. From Player A's perspective, though, the decision epochs correspond to states $s$ and $s'$, not $\Bar{s}$. In our numerical results we defined $N=1000$ and considered integer values of $\epsilon$ from 1 to 20, resulting in over 75 million simulated shots.

\vspace{0.5cm}

% \begin{center}
\begin{algorithm}[H] 

    \For{\emph{state} $s \in S$}
    {
        Generate $N$ shots $w_1, \ldots, w_N \sim G \: | \: s$
        
        Fit the intention distribution $f_s$
        
        For each $a \in A$, fit an execution distribution $\mathcal{N}(\bm{\mu}_{s,a}, \bm{\Sigma}_{s,a})$ %for all $a \in A$ \tcyc{this needs to be tweaked...it depends on $\epsilon$ but there is no $\epsilon$ defined at this point yet} \cf{maybe this $N_2$ doesn't need an epsilon in it, since its essentially $\epsilon=1$ here?}
        
        \For{$\epsilon = 1, 2, \dots, E$}
        {
            \For{$i = 1, 2, \ldots, N$}
            {
                Draw an intention $a \sim f_s$
                    
                Draw a landing location $(x, y) \sim \mathcal{N}(\bm{\mu}_{s,a},\epsilon\bm{\Sigma}_{s,a})$
                
                Generate a shot $w \sim G \:|\: s,(x,y)$ and record $\Bar{s}$
            
                Generate a return shot $w' \sim G \:|\: \Bar{s}$ and record $s'$
            }

            Compute $P_{\epsilon}(s' \: | \: s,a)$
        }
    }
\caption{Generating Transition Probabilities with Error}
\label{algo:TransProbs}
\end{algorithm}

\vspace{0.5cm}

% \end{center}

\textbf{Validating ``average'' error:} To determine which value of $\epsilon$ represents the error level of an ``average'' player from our simulated shot data, we examine the empirical (i.e., the VON CRAMM-generated shots without any fitting or scaling) probability of the three outcomes of a shot -- win the point, lose the point, or continue (i.e., remain in play) -- and compare it with their simulated values for different $\epsilon$ drawn from the fitted distributions. Intuitively, the probability of winning the point outright or continuing should decrease with $\epsilon$ and the probability of losing the point outright should increase with $\epsilon$. This intuition is reflected in Table \ref{table:averageError} in Appendix \ref{sec:AppendixA}. Given that the probability of losing the point and keeping the point in play were monotonically increasing and decreasing, respectively, we chose the value of $\epsilon$ that generated probabilities close as possible to the empirical probabilities ``without going over'', which turned out to be $\epsilon = 13$. Choosing this value also avoids biasing towards out-of-bound shots.  To be clear, this is a different scaling of the covariance matrix than what was implemented when we fit the perfect execution distribution. The perfect execution scenario corresponds to $\epsilon = 1$, which corresponds to the covariance matrix where 90\% of the shots land in a given intention. Further note that each intention has its own covariance matrix that was scaled differently to define perfect execution. But once those covariances are defined for perfect execution for each intention, then this single $\epsilon$ scaling is used to model growing execution error over all intentions.

% \rs{[I am missing something here. Didn't you calibrate the Sigma matrix using VON CRAMM? If so, shouldn't this exercise give you epsilon = 1?  Is the scaling down to 90\% perfection resulting in epsilon = 13?]} \cf{[Same nuance as above, idk if we need to be more specific or if a particular wording we're using is confusing]} 

\begin{figure}[H]
\centering
\begin{tabular}{cccc}
    \includegraphics[width=12.5cm]{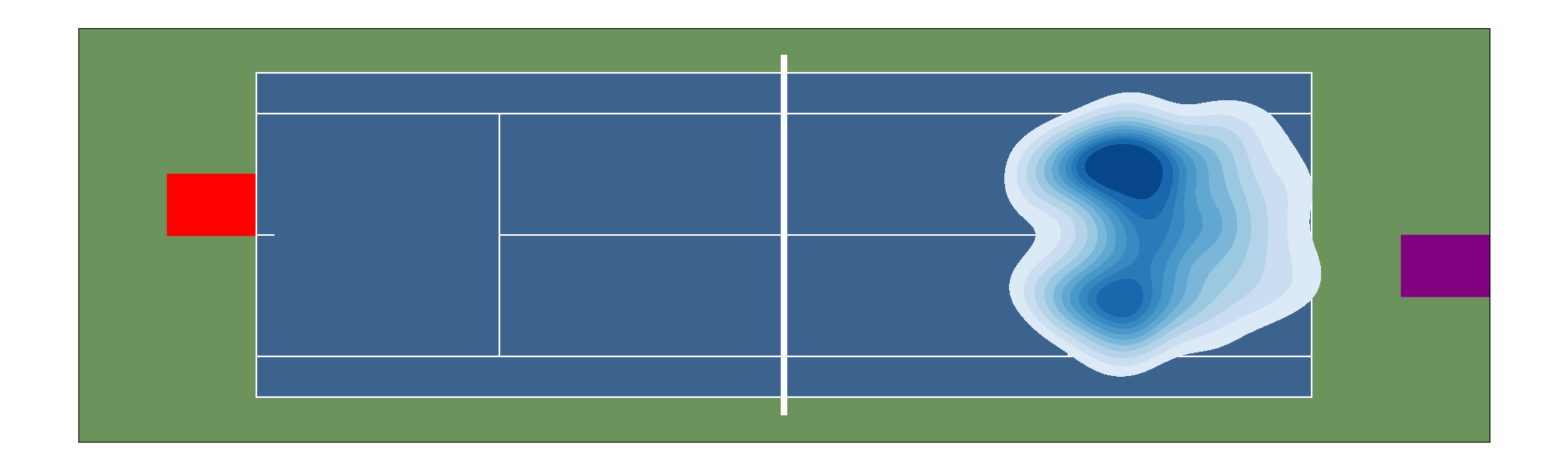} \\
    (a) \\[6pt]
\end{tabular}
\begin{tabular}{cccc}
    \includegraphics[width=12.5cm]{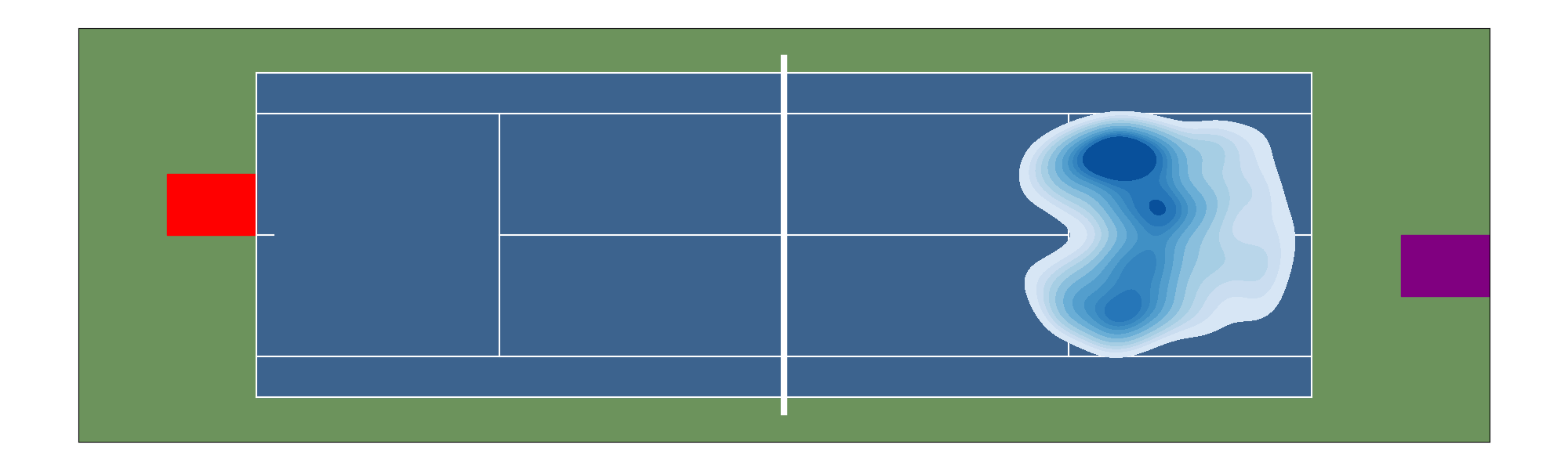} \\
    (b) \\[6pt]
\end{tabular}
\begin{tabular}{cccc}
    \includegraphics[width=12.5cm]{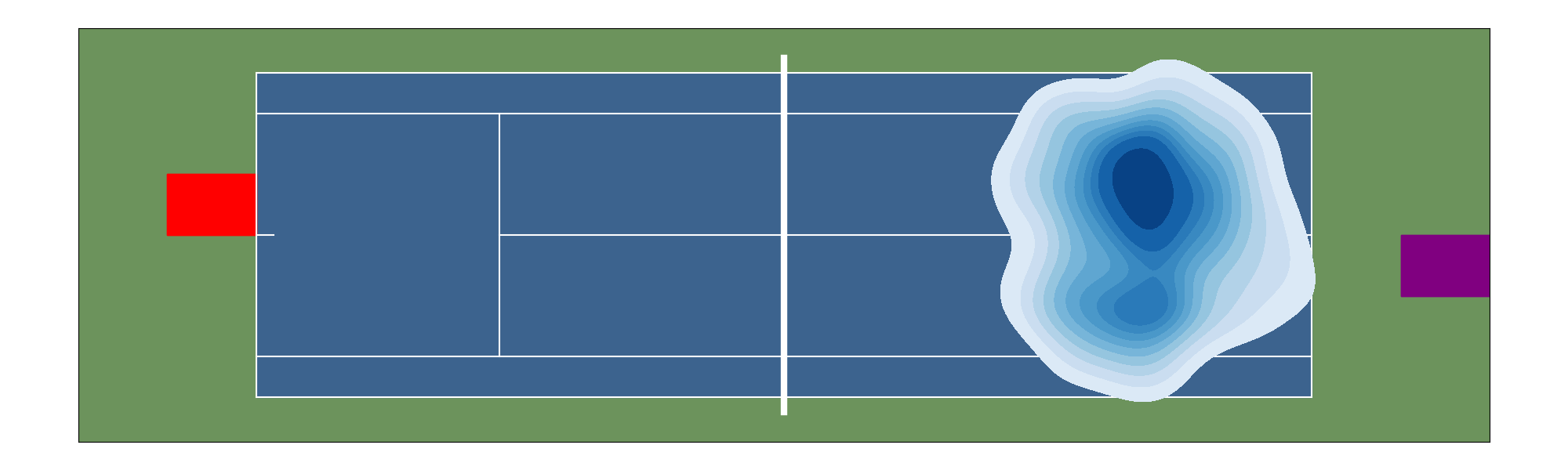} \\
    (c) \\[6pt]
\end{tabular}
\caption{(a) The raw landing location distribution over all actions when Player A is in cell 10 (red square), Player B is in cell 81 (purple square) and Player A is hitting a rally shot. (b) The simulated landing locations for $\epsilon=1$. (c) The simulated landing locations for $\epsilon=13$, the closest fit to the empirical average player.}
\label{fig:errormodeling}
\end{figure}

To illustrate, consider the state $s = (10,81, \text{rally})$. Figure \ref{fig:errormodeling} shows the empirical shot distribution, the simulated shot distribution for $\epsilon=1$, and the simulated shot distribution for $\epsilon=13$ over all actions. We can see that the $\epsilon=1$ case models a player who almost never hits it out of bound, representing the perfect execution case. For $\epsilon = 13$, although the simulated distribution does not perfectly match the empirical distribution, it mostly recovers the salient features, namely the existence of two modes, the location of the primary mode, and the general shape of the distribution's support.
%, there is a strong alignment between the empirical results and the simulated results for $\epsilon=13$, further validating the appropriateness of deeming $\epsilon=13$ to model the ``average" player. 
Note that while we can tune the execution error of Player A through $\epsilon$, Player B's execution error remains fixed at the ``average'' error level. In other words, when we choose $\epsilon = 13$, we are modeling a point between two symmetric players.
%an implicit consequence in our approach is that Player B remains at a constant average error level. Only Player A's execution error varies with the choice of $\epsilon$. Therefore, when $\epsilon=13$, we are able to model two players, both of average error-level, facing one another.

% \tcyc{looking at the table, it seems eps 14 could be argued to be better than eps 13} \cf{I chose eps=13 as it has lower P(error) than eps=14. This is to avoid biasing out-of-bound shots}
%we compare the three main types of outcomes for a shot: win the point, lose the point, remain in-play. Computing the probabilities for these outcomes for each $\epsilon$ and comparing them to the empirical probabilities of average players, we see a strong relationship at $\epsilon=13$ (see 

% \tcyc{the pictures don't look all that similar to me between empirical and eps = 13. can we perhaps find a different state where the qualitative images look more similar between these two cases? also, the out of bounds prob seems noticably higher for eps 13 than empirical. } \cf{This is the closest fit I could find after 60-70 attempts. I think whats nice about this one is the bi-modal nature is captured. There is more OB, but I also liked that in both cases the corners are untouched. Maybe more text highlighting that the fits aren't perfect on the edges, but there is greater alignment in places of high concentration?} 

\subsubsection{Rewards} \hfill\\ \vspace{-0.5cm}

Since the point is modeled from the perspective of Player A, reaching the winning state for the first time results in a reward of $+1$, reaching the losing state results in a reward of $0$, and the reward for all other states is 0. That is, $r(s,a,s') = 1$ when $s \ne W$, $s' = W$ and $r(s,a,s') = 0$ otherwise.

\subsection{Bellman equation}

We now discuss how these attributes are used in the MRP and MDP models, as well as how each model is solved.

\subsubsection{Markov reward process}\hfill\\ \vspace{-0.5cm}

In an MRP, the goal is policy evaluation. That is, we assume a given policy and we compute the value function for this policy using the Bellman equation as a system of linear equations. We will use the MRP model to compute value functions for different types and magnitudes of execution error in order to study how these errors impact a player's probability of winning the point. As such, we will solve the MRP with different transition probability matrices, each corresponding to a different level of execution error. The particular policy we focus on is the randomized policy corresponding to the intention distribution $f_s$. %-- the intention distribution of the simulated shot locations from $G$. 
%\tcyc{want to look later to make sure all uses of the word `empirical' are consistent.} 
This policy is randomized because given any state $s$, there is a distribution of shot landing locations over the action set $A$, which means a nondegenerate probability distribution describes the possible actions chosen in state $s$. Let $\hat\pi$ denote this policy. Then, the value of this policy for execution error level $\epsilon$ can be computed according to %\cf{I fixed equations below, the n-step ahead may need to be tweaked to capture the $n$ aspect.}
\begin{align} \label{eq:MRP-value}
    % V_\epsilon^{\pi^e}(s) & = \mathbf{E}_{f_s}\left(\mathbf{E}_{P_\epsilon}\left(r(s,a,s') + V_\epsilon^{\pi^e}(s')\right) \right), \quad \forall s \in S \\
    %& = r(s,a,s') + \sum_{a \in A}  P(a \: | \: s)\left(\sum_{s' \in S} P_{\epsilon}(s' \: | \: s,a)\cdot V_\epsilon^{\pi^e}(s') \right), \qquad \forall s \in S  \\
    V_\epsilon^{\hat\pi}(s) & = \sum_{a \in A} f_s(a)\left( \sum_{s' \in S} P_{\epsilon}(s' \: | \: s,a)\left(r(s,a,s') +  V_\epsilon^{\hat\pi}(s')\right) \right), \quad \forall s \in S 
\end{align}
where $V_\epsilon^{\pi}(s) = 0$ when $s$ is an absorbing state ($W$ or $L$) for any policy $\pi$. 

% %We first explore the MRP where we have a probabilistic action distribution for each state. This can be seen as the value for the ``empirical'' policy used by players, and is not subject to optimization. 
% We introduce a new probability, $P(a \: | \: s)$, as the probability of choosing action $a$ in state $s$. \tcyc{this is the empirical intention distribution for state $s$, right? if so, we should define this earlier and write the math in Alg 1 and in prev sections} \cf{yeah, this would be $f_s$} These probabilities can be understood as the empirical proportions that is seen in the data, and as such, we define this policy as $\pi^e$. This is modelled by the Bellman equation shown in \ref{eq:MRP-value}. Note that we define $\Bar{S} = S \setminus {W}$ to indicate all states except the state indicating Player A won the point, as this is related to the reward of $+1$. 
% \begin{equation} \label{eq:MRP-value}
%     V^{\pi^e}(s)= \sum_{a \in A}  P(a \: | \: s)\left(p_{s,W}^a + \sum_{s' \in \Bar{S}} p_{s,s'}^a \cdot V^{\pi^e}(s') \right) \qquad \forall s \in S 
% \end{equation}
% \noindent
% The $V$'s in this case represent the probability Player A wins the point when starting in state $s$ and following the empirical policy $\pi^e$. \\

\subsubsection{Markov decision process} \hfill\\ \vspace{-0.5cm}

In an MDP, the goal is policy optimization. We use the MDP model to investigate how the presence of execution error affects the structure of the optimal policy and the optimal value function. An optimal value function for execution error level $\epsilon$ can be computed through the following optimality equation
\begin{align}
    \label{eq:MDP-value}
  %  V_\epsilon^{*}(s) & =\max_{a \in A}\left\{\mathbf{E}_{P_\epsilon}\left(r(s,a,s') + V_\epsilon^{*}(s')\right) \right\}, \quad \forall s \in S \\
    V_\epsilon^{*}(s) &= \max_{a \in A} \left\{\sum_{s' \in S} P_{\epsilon}(s' \: | \: s,a)\left(r(s,a,s') +  V_\epsilon^{*}(s')\right)\right\} , \quad \forall s \in S. 
\end{align}
Let the corresponding optimal deterministic policy be denoted $\pi_\epsilon^*$, which is defined by the following equation
\begin{equation}
    \label{eq:opt-policy}
    \pi_\epsilon^*(s) \in \argmax_{a\in A}  \left\{\sum_{s' \in S} P_{\epsilon}(s' \: | \: s,a)\left(r(s,a,s') +  V_\epsilon^{*}(s')\right)\right\}, \quad \forall s \in S. 
\end{equation}

In practical terms, the MDP models the situation in which each state has one optimal action that is always chosen when Player A enters that state. For example, if the optimal action in state $(8,75, \text{rally})$ is to aim for the deuce corner (CD1), then the MDP will model a game where every time Player A reaches $(8,75, \text{rally})$, they will aim for CD1. This deterministic policy does not appropriately represent the nature of tennis as Player A will become fully predictable. Moreover, it may not be possible to choose an optimal action on every shot due to inherent shot difficulty, which is what is implied if the policy $\pi_{\epsilon}^*$ is used. Thus, the value computed from \eqref{eq:MDP-value} may be unrealistic and overly optimistic.

Therefore, we consider an intermediate approach where we implement a series of optimal actions determined from $n$ one-period models (i.e., optimizing only for one decision epoch) and then ``roll-out'' with policy $\hat\pi$. In essence, this is policy iteration starting with the policy $\hat\pi$ but terminating after $n$ iterations.
%
%Therefore, we also consider a $n$-step lookahead policy with respect to the value function $V_{\epsilon}^{\pi^e}$. \cf{[whatever we call this new model, we need to update its name throughout the paper.]} While such an approach is typically used as a way to more tractably, but approximately, solve an MDP, our purpose is slightly different. 

We consider such an approach to explicitly model a player choosing a sequence of greedy (i.e., one-period optimal, assuming a fixed value function for the cost-to-go) actions, before reverting to the policy corresponding to the intention distribution $f_s$ for the remainder of the shots in the point. %\rs{[Why do you call it "greedy"? In (4), the long-run value is captured (but using the baseline policy $\pi^e$). I feel "greedy" might confuse the reader to believe that the long-run value is ignored.  Same comment with the usage of "one-period models" above.]} \cf{[I think when we use greedy or one-period, we mean we are just \emph{optimizing} for one period. Is that fair to say?]} 
This setup allows us to isolate the value of optimal action selection for a limited number of shots, which provides a more realistic comparison with the the impact of varying execution error. It also provides insight into which shots in the point provides the largest incremental gain, if an optimal action is selected, and the magnitude of this gain. 

%\tcyc{need to think about this more, if we just position this as policy iteration iterations instead of n-step}

We start with $n = 1$ and the corresponding one-period greedy policy, denoted $\pi_\epsilon^1$. This policy can be written as $\pi^1_\epsilon = (d^1, \hat d, \hat d, \ldots)$ where $\hat d(s)$ is the decision rule that chooses an action according to the intention distribution in state $s$ and $d^1$ is determined from

\begin{equation}
    \label{eq:MDP-1step:argmax}
    d^1(s) \in  \arg\max_{a \in A} \left\{\sum_{s' \in S} P_{\epsilon}(s' \: | \: s,a)\left(r(s,a,s') +  V_\epsilon^{\hat\pi}(s')\right)\right\} , \quad \forall s \in S .
\end{equation}
The value of this policy under execution error level $\epsilon$ is computed according to %\cf{[For the mdp policy and the n-step policy on the decision epochs, we need to denote the eps as well. I tried adding that]}

%\cf{We note that in practice, players would not choose an optimal action on every shot they try to hit. This is because players often use seemingly sub-optimal shots that are intended to ``set up a winning shot".} Therefore, we will also implement a $n$-period lookahead version of the MDP model, not as a way to approximately solve the main MDP model, but as a way to isolate the impact of having perfect execution for one shot or two shots, etc. \cf{This allows us to model a player being selectively optimal and provides insights into which moment in the point provides the largest incremental gain, if played optimally.} This $n$-period lookahead model seeks to optimize the current period and then follow the empirical policy ($\pi^e$) for the remainder of the point (as in the MRP model). It can be computed through the following optimality equation

\begin{equation}
    \label{eq:MDP-1step}
    V_\epsilon^{\pi_\epsilon^1}(s)= \max_{a \in A} \left\{\sum_{s' \in S} P_{\epsilon}(s' \: | \: s,a)\left(r(s,a,s') +  V_\epsilon^{\hat\pi}(s')\right)\right\} , \quad \forall s \in S. 
\end{equation}
Note that if $\hat\pi$ is replaced with $\pi_\epsilon^*$, then $\pi_\epsilon^1 = \pi_\epsilon^*$. In general, for $n \ge 1$, we compute $d^n$ from

\begin{equation}
    \label{eq:MDP-nstep:argmax}
    d^n(s) \in \arg\max_{a \in A} \left\{\sum_{s' \in S} P_{\epsilon}(s' \: | \: s,a)\left(r(s,a,s') +  V_\epsilon^{\pi^{n-1}_\epsilon}(s')\right)\right\} , \quad \forall s \in S, 
\end{equation}
where $\pi^{n-1}_\epsilon = (d^1, d^2, \ldots, d^{n-1}, \hat d, \hat d, \ldots)$. The value of the policy $\pi^n_\epsilon = (d^1, d^2, \ldots, d^{n}, \hat d, \hat d, \ldots)$, composed of $n$ greedy actions followed by the randomized policy corresponding to the intention distribution $f_s$, is computed via 
\begin{equation}
    \label{eq:MDP-nstep}
    V_\epsilon^{\pi_\epsilon^n}(s)= \max_{a \in A} \left\{\sum_{s' \in S} P_{\epsilon}(s' \: | \: s,a)\left(r(s,a,s') +  V_\epsilon^{\pi^{n-1}_\epsilon}(s')\right)\right\} , \quad \forall s \in S 
\end{equation}
with the boundary condition that $\pi^0 := \hat\pi$. It follows that for any positive integer $n$, the value function associated with $\pi_\epsilon^n$ is bounded between the MRP value function and the MDP value function, i.e., $V^{\hat\pi}(s) \le V^{\pi_\epsilon^n}(s) \le V^*(s)$ for all $s \in S$.

\subsubsection{Value function interpretation.}

For any policy $\pi$, the value function $V_\epsilon^{\pi}(s)$ has a natural interpretation within the context of tennis.

\begin{thm}\label{thm:V_mrp}
	Given any policy $\pi$, $V_\epsilon^{\pi}(s)$ equals the probability that Player A wins the point given that the current state is $s$ and under execution error level $\epsilon$.
\end{thm}

The proof of this result is straightforward and omitted. In essence, under a fixed policy $\pi$, the induced stochastic process becomes a Markov chain, and given the reward function structure, $V^\pi(s)$ becomes the probability of absorption to the win state. See \citet{Chan2016} for more details. Given this result, subsequent references to the ``value'' of a policy should be interpreted as the probability of winning the point starting from the given state. 

%\tcyc{maybe here we want to introduce the concept of V(start)?} \cf{I added it in sec5.2, I think it fits better there.}

% \proof{}
% 	Given a policy $\pi$, the Bellman equations for $V_\epsilon^{\pi}(s)$ \eqref{eq:MRP-value} are exactly the system of equations describing the probability of absorption in a Markov Chain \citep{Bertsekas2008} via the win state. Therefore, $V^{\pi}(s)$ represents the corresponding probability under policy $\pi$. \Halmos
% \endproof

% \begin{cor}
% \label{cor:V_MDP} $V_\epsilon^*(s)$ equals the probability that Player A wins the point when in state $s$ and using the optimal policy $\pi^*$ under execution error level $\epsilon$.
% \end{cor}

\section{Results}

In this section, we present results from several numerical experiments applying our MRP and MDP frameworks. In particular, we 

%\vspace{-0.75cm}
\begin{enumerate}
    \item Quantify the loss in value as $\epsilon$ increases (Section \ref{sec:MRP_Value}),
    \item Explore how, and to what extent, particular types of error affect a player's value (Section \ref{sec:MRP_Value_and_Error}),
    \item Compare how aggressive versus conservative players perform given differing magnitudes of execution error (Section \ref{sec:MRP_Playstyle}),
    \item Quantify the value a player can achieve by choosing optimal actions on all shots, and how optimal shot selection changes in the presence of execution error (Section \ref{sec:MDP_Results}),
    \item Determine how the value changes when a player can choose optimal actions on only a subset of shots (Section \ref{sec:MDP-nStep}).
\end{enumerate}

%\tcyc{if we can describe each in one sentence, then do the numbered list with (Sec xxx) at the end of each one. If we can't, then write out a paragraph where each analysis/subsection is one or two sentences on their own} \cf{Done above - just a note: sections 51. and 5.4 don't have experiments per se but is rather just displaying the MRP/MDP results.}

\subsection{Relationship between value and execution error}
\label{sec:MRP_Value}
\noindent 

\begin{figure}[H]
\centering
\begin{tabular}{cccc}
    \includegraphics[width=0.48\textwidth, height=4cm]{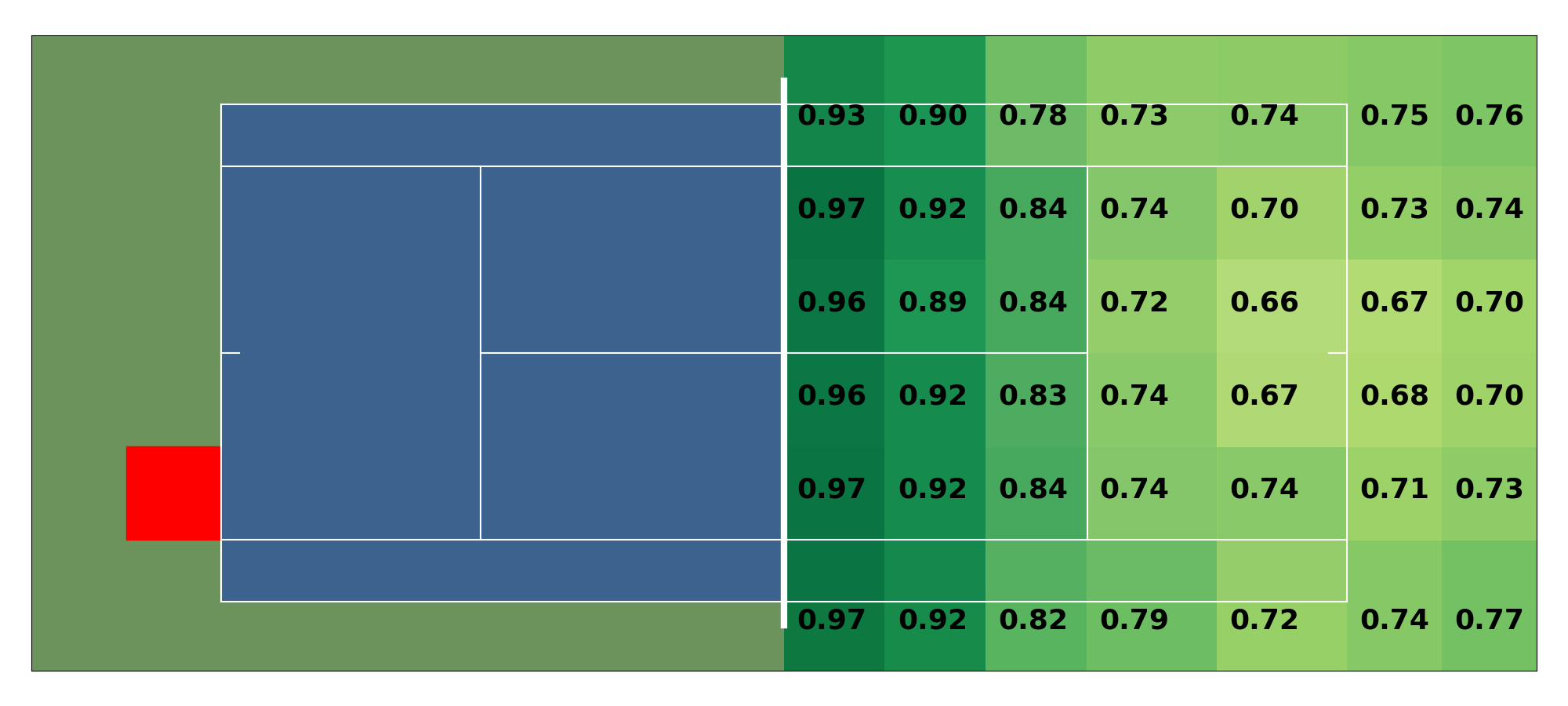} &
    \includegraphics[width=0.48\textwidth, height=4cm]{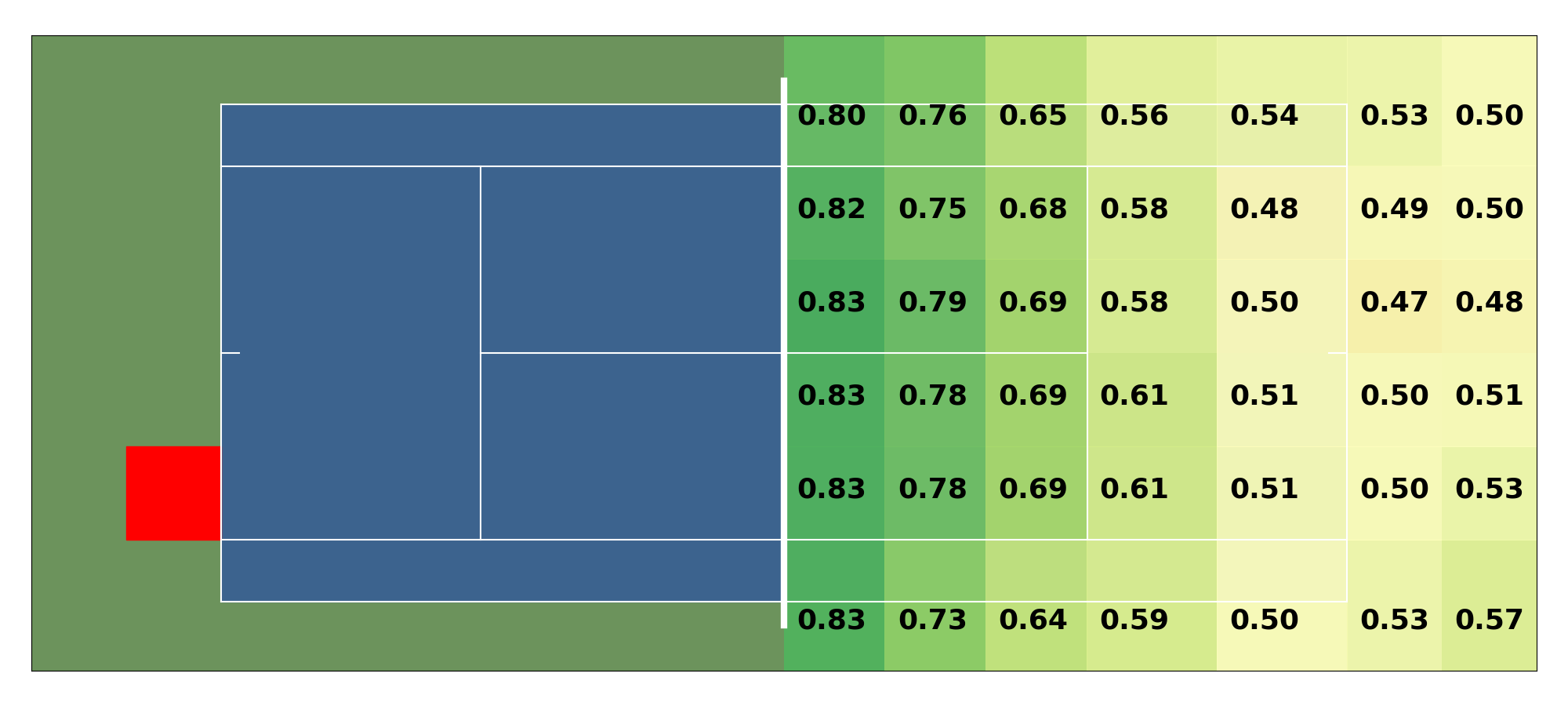} \\
    (a) Player A at deuce baseline with $\epsilon = 1$ & (b) Player A at deuce baseline with $\epsilon = 13$\\[6pt]
\end{tabular}
\begin{tabular}{cccc}
    \includegraphics[width=0.48\textwidth, height=4cm]{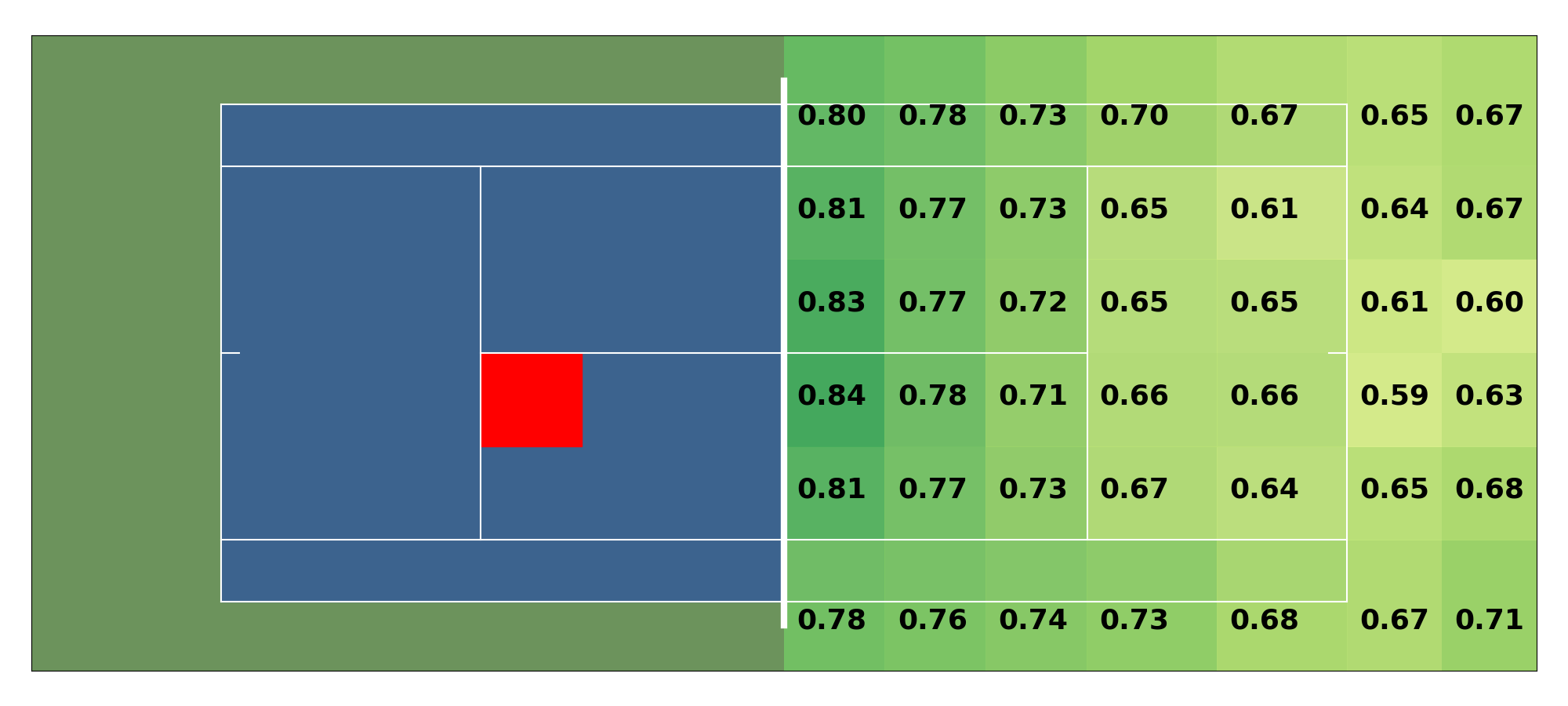} &
    \includegraphics[width=0.48\textwidth, height=4cm]{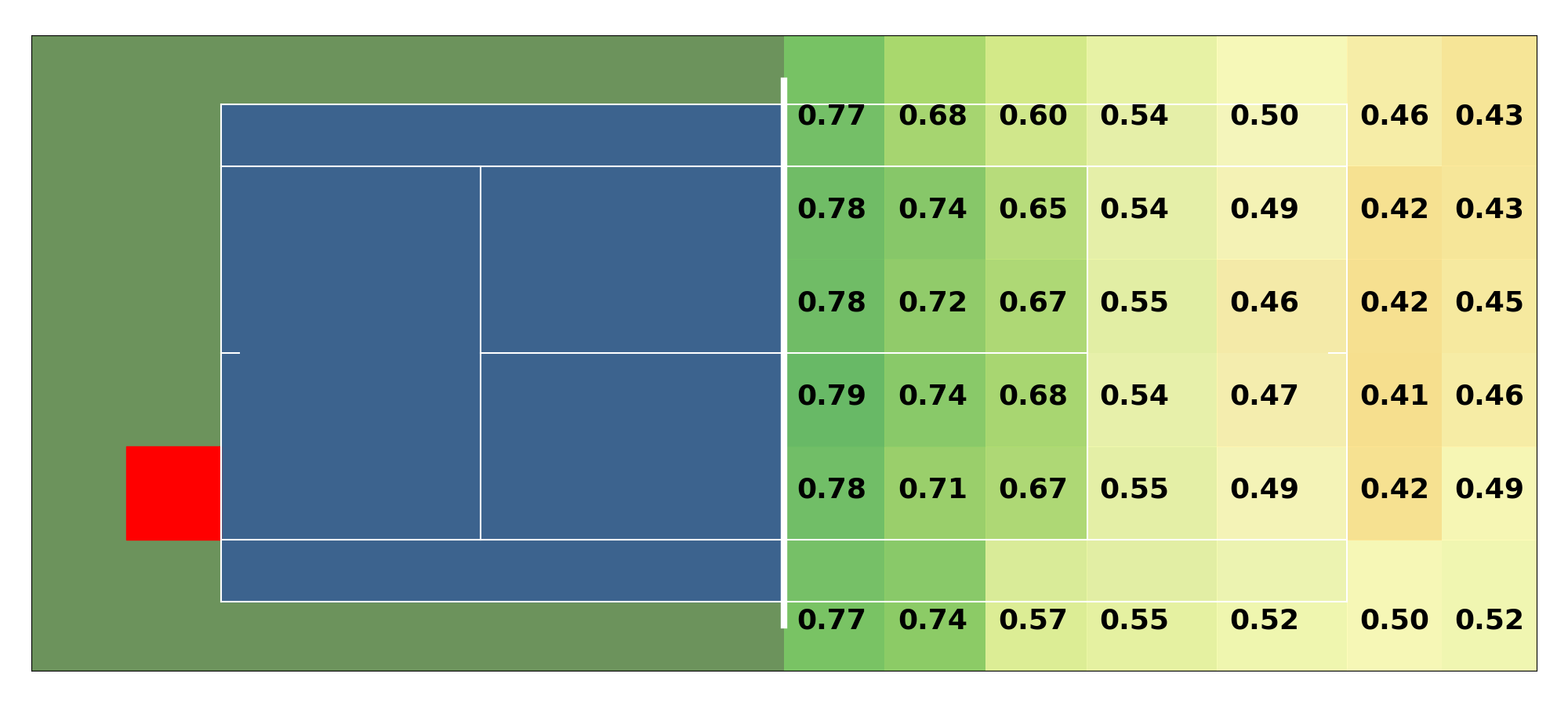} \\
    (c) Player A at deuce service line with $\epsilon = 13$ & (d) Player A at deuce baseline with $\epsilon = 20$\\[6pt]
\end{tabular}
\caption{The MRP value function for various rally states and values of $\epsilon$. The number in each cell on the right side of the court indicates the value of that state, where the state is described as Player B standing in that cell, Player A standing in the red cell, and the shot being a rally shot being hit by Player A. %\tcyc{when you're adjusting the font/numbers in the other fig, maybe make these a bit bigger too? get rid of the eps at the top of the fig. we will put it in the caption}
}
\label{fig:mrp_Vs}
\end{figure}

Figure \ref{fig:mrp_Vs} illustrates the MRP value function for several different rally shot states and for different values of $\epsilon$. The number in each cell on the right side of the court indicates the value of that state, where the state is described as Player B standing in that cell, Player A standing in the red cell, and the shot being a rally shot being hit by Player A. For example, in Figure \ref{fig:mrp_Vs}(a), when Player B is standing near the opposite baseline corner as Player A, there is a roughly 70\% chance that Player A will win the point. This value is much greater than 50\% because it corresponds to the situation where Player A can execute their shot perfectly ($\epsilon = 1$). In Figure \ref{fig:mrp_Vs}(b), with average execution error ($\epsilon=13$), the value drops to roughly 50\%. As a reminder, the execution error of Player B remains fixed at the average level of $\epsilon = 13$ throughout. Thus, this result is intuitive since it models two equally skilled players trading cross-court rally shots, each having an equal chance of winning the point. When Player A gains net control in Figure \ref{fig:mrp_Vs}(c), the values appropriately increase compared to Figure \ref{fig:mrp_Vs}(b). Lastly, going back to the cross-court rally scenario in Figure \ref{fig:mrp_Vs}(d), the value function drops below 50\% since Player A has higher execution error than average ($\epsilon = 20$). 

\subsection{Impact of particular types of execution error}
\label{sec:MRP_Value_and_Error}

Next, we drill down into how a player's error impacts the value function. To quantify the effect of a particular error scenario, defined as a specific error type (i.e., the particular states in which player's experience execution error) and value of $\epsilon$, we build transition probability matrices that only model $\epsilon$-level error in those states. For example, to isolate the effect of a player experiencing an execution error of $\epsilon = 7$ on ad serves, we modify the transition probability matrix so that it has transition probabilities corresponding to $\epsilon=7$ for the ad serve states and transition probabilities of $\epsilon=1$ for all other states. We will use this approach to model the situation where, starting from perfect execution, we grow the execution error towards average error and then beyond in a specific state. We will also grow the execution error in all states simultaneously, which we refer to as the ``total error'' scenario.  

%\tcyc{or should it be eps = 13 on all the other states? i guess we just have to explain clearly this is error is wrt perfect execution or wrt average execution. } \cf{I think eps=1 makes the most sense in terms of the graphs, as all the lines will start with perfect execution and grow towards average and then beyond average error. Otherwise we'll have the case where say we're analyzing deuce serves, players are average everywhere but perfect on deuce serves, then growing deuce serve error to the average and then beyond average. eps=1 seems to be the cleanest to me} This is repeated for all levels of $\epsilon$ and all the specific scenarios to be tested. We also include the case where execution error grows in every state, denoted as the ``total error" scenario. \tcyc{at this point the prev two sentences are not clear enough to be understood w/o reading ahead} \cf{My goal for this paragraph was just outline the different scenarios we'll be testing - I don't think the metric of comparison is necessary at this stage?}

% For each scenario, defined as a specific error type and value of $\epsilon$, \tcyc{i like this definition of a scenario, but couldn't quite fit it in below. will revisit later}

To make the impact of execution error comparable across different states (e.g., rally states vs. serve states), we roll back the state value to the start of the point, since the value of subsequent states is captured through the Bellman equation \eqref{eq:MRP-value}. For example, if increased execution error causes the value of a rally state to decrease, this will be reflected in a decrease in value on serve states. To determine an overall value for all serve states, we compute a weighted average based on the empirical frequency of points starting in each serve state. This approach is repeated for the serve return states. Since Player A will never transition between a serve to a serve return in the same point, these two states are independently impacted by the value on the rally states. Assuming that Player A is equally likely to start the point as the server or serve returner, we average the two values for the serve states and the serve return states to compute the catch-all value, which we denote $V(\text{starting-state)}$. It is important to note that Player A's first shot in the point will always be either the serve or serve return, and their second shot will always be a rally.

First, as a rough validation, note that the weighted average value on serves for the total error scenario with $\epsilon = 13$ is roughly 0.652. This can be interpreted as the probability a player with average execution error wins the point when serving, according to our model. This result is in strong agreement with an empirical point-win probability on serve, 0.642, computed from a large dataset of the Australian Open from 2017-2019 \citep{Sackman2021}. %\tcyc{(comprising xxx data points from Grand Slam matches. Is it men's or womens?)} \cf{addressed in previous sentence.} 

\begin{figure}[H]
    \centering
    \includegraphics[width=\textwidth]{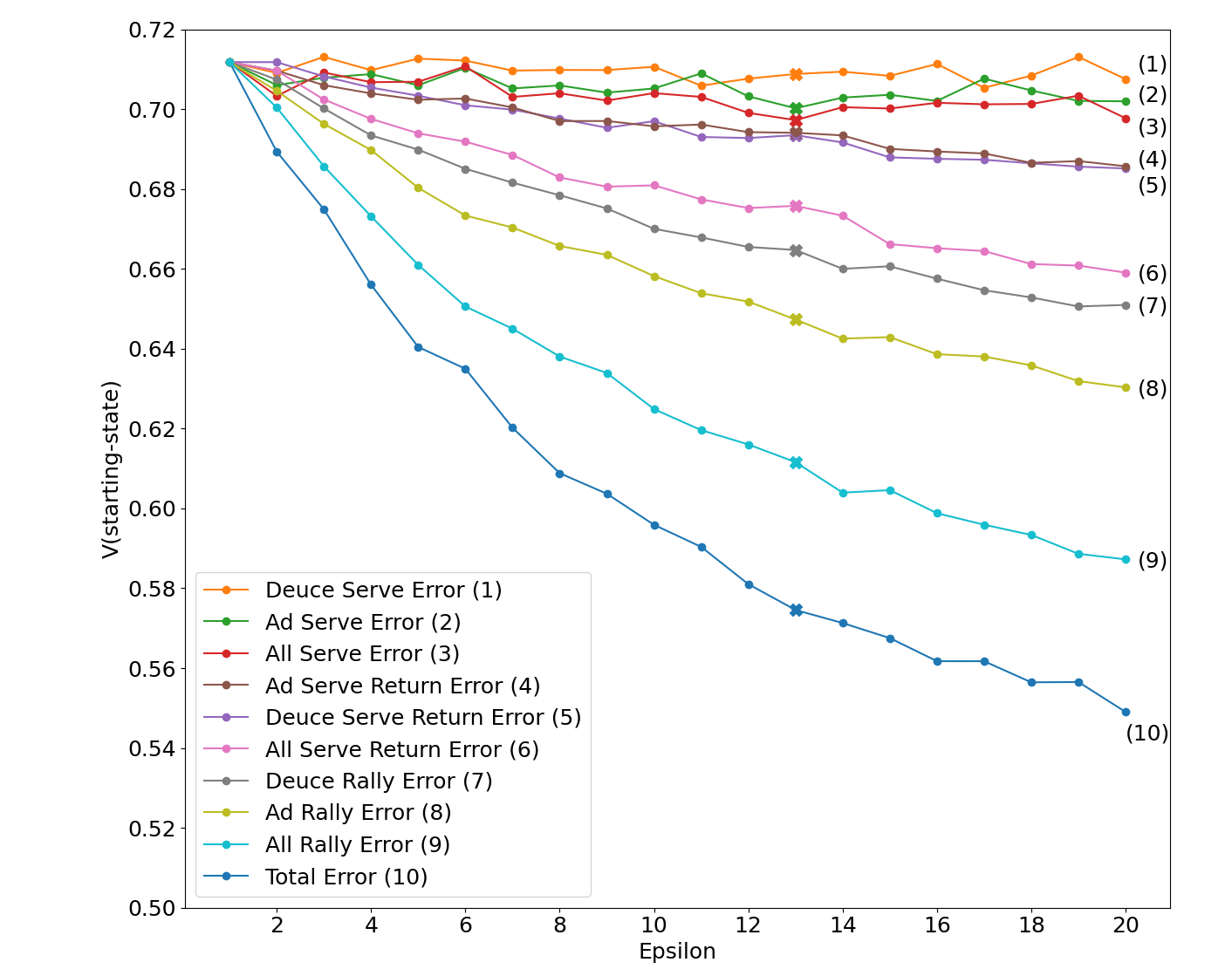}
    \caption{The relationship between value and error for various scenarios. Value is measured as the weighted average for the value on both the serve and serve return states, denoted, $V(\text{starting-state)}$. Each scenario depicts the value when Player A has perfect execution in all states, and epsilon-level execution error in that particular scenario.}
    \label{fig:mrp_errors}
\end{figure}

Figure \ref{fig:mrp_errors} shows the relationship between $V(\text{starting-state)}$ and different error scenarios. Note, the non-monotone nature of some of the curves is due to finite sampling. %\tcyc{[did we say size of simulation? i also hesitate say anything about computational expense for error bars, since if we really wanted to do it, we'd just be running it now (or in review), and maybe not for all epsilon for all lines, but just a few to give an idea.]} \cf{[We mention the 75 million shots on the bottom of page 8. We never mention the time required for this though. I agree, I think your shorter comment on the non-monotone behaviour is sufficient.]} % Unfortunately, the computational expense required to generate error bars from the simulation jagged behaviour in the lines is due to sampling error. The computation expense necessary to provide error bars for a simulation of this size was deemed infeasible.} \rs{[Thoughts on why some curves are non-monotone?]} \cf{[tried to mention the sampling error above]}
%
% \tcyc{we should say a few other high level comments about the fig. describe what it is showing. how to read it. what the y-axis shows (eg are there two weighted avgs in c)?) . xaxis is self explanatory. comment on serve being higher value than serve return in general? also for all states, it is between serve and serve return. but some of the states in the all states plot only belong to one of the two. eg "all serve return error" only shows up in serve return states, and not serve states. i'm trying to understand how the pink line in all states is higher than the pink line in serve return states. there is some averaging wrt serve states, but how is it done?} \cf{just put (C)}
%
One insight is that execution error in ad rallies appear to be more detrimental than in deuce rallies. We define an ad (deuce) rally shot as a shot that originated from the ad (deuce) side of the court, regardless of the proximity to the net. To understand this finding, we first note that the VON CRAMM framework models performance aggregated over all players. There is roughly a 90\% prevalence of right-handed players in tennis. As such, it is reasonable to assume that in the majority of cases, an ad (deuce) rally shot represents a backhand (forehand) shot. Moreover, across all error scenarios, given that the point did not end in an immediate win or loss, 55\% of the time it transitioned into an ad rally state compared to deuce rally state. This means that the outcome of shots from an ad rally state will have a larger impact on $V(\text{starting-state)}$. Looking at the transition probabilities from ad and deuce rally states to the absorbing state, they are indistinguishable except for transitions into the state where Player A errors in the shot (see Appendix \ref{sec:AppendixB}). Essentially, as $\epsilon$ grows, the chance of Player A hitting the ball out of bounds grows faster for ad rally shots than deuce rally shots. %Considering that backhanded shots are harder to hit than forehands, and the prevalence of right-handed players, this finding is compatible with tennis intuition. 
While there appears to be a difference between ad and deuce rally shots, this discrepancy in value does not seem to exist when comparing ad and deuce serves or serve returns.

A second finding is that execution error on serves does not have a large impact on overall value. While stellar serves can lead aces or more difficult returns for the opponent, this likely has as much to do with speed and spin, as it does with a tight execution error around landing location. As such, growing the $\epsilon$ value in our modeling framework for serves doesn't impact value substantially. If we expanded our definition of execution error to include deviations from a player's intended spin, speed, and bounce angle, in addition to the intended landing location, we would expect to see further loss in value as execution error increased. Additionally, if we modeled a finer action set on serves that captured the corners of the service box, for example, then reducing execution error on these actions would likely improve value. These extensions are beyond the scope of this paper and we leave it as future work. %\rs{[Does this finding remain the same if action set for serve is made finer? (Referring to Fig 2b.) As a tennis player, I am not convinced that this is true. Location is definitely important (as well as speed and spin). Pros are able to return very fast balls; it is when the ball is placed at the corners (and is served fast) that results in an ace.]} \cf{[tried to fix above]}

%Lastly, execution error on rally shots destroy more value than execution error on serve return shots.  since rally shots can happen repeatedly in a point, whereas serve returns are just once, compounding poor execution on rallies can greatly impact overall value. \cf{[Having trouble explaining this...]}

% \begin{figure}[H]
% \centering
% \begin{tabular}{cccc}
%     \includegraphics[width=0.30\textwidth, height=7.5cm]{Figures/mrp_error_serve.png} &
%     \includegraphics[width=0.30\textwidth, height=7.5cm]{Figures/mrp_error_serve_return.png} &
%     \includegraphics[width=0.30\textwidth, height=7.5cm]{Figures/mrp_error_starting_state.png} \\
%     (a) Serve states & (b) Serve return states & (c) All states \\[6pt]
% \end{tabular}
% \caption{The relationship between value and error for various scenarios. Value is measured in three ways: the weighted average for the value on the serve states, the serve return states, and both serves and serve returns. Each scenario depicts the value when Player A has perfect execution in all states, and epsilon-level execution error in that particular scenario. \tcyc{can we order the legend so it lines up with the ordering of the lines in the figs?}}
% \label{fig:mrp_errors}
% \end{figure}

\subsection{Impact of execution error on play styles}
\label{sec:MRP_Playstyle}
%\tcyc{PAUSED HERE}

All of the previous analysis thus far assumed a player with a constant intention distribution but varying magnitudes and types of error. We now explore how execution error may affect players with different play styles. First, we define a ``conservative'' shot as one that is aimed at an intention in the middle of the court (any intention starting with `M'). Any other shot (short, sideline, baseline and corner) will be considered ``aggressive''. Then, we define a conservative (aggressive) player as one who hits 50\% more conservative (aggressive) shots than the player analyzed in previous sections (i.e., a player who selects shots in accordance with the intention distribution $f_s$). We can model this change in play style by adjusting the intention distributions in Algorithm \ref{algo:TransProbs}. While players are likely to be more aggressive or conservative in only certain states, for simplicity we examine the situation where they are uniformly more aggressive or conservative in all states. 

\begin{figure}[H]
    \centering
    \includegraphics[width=\textwidth]{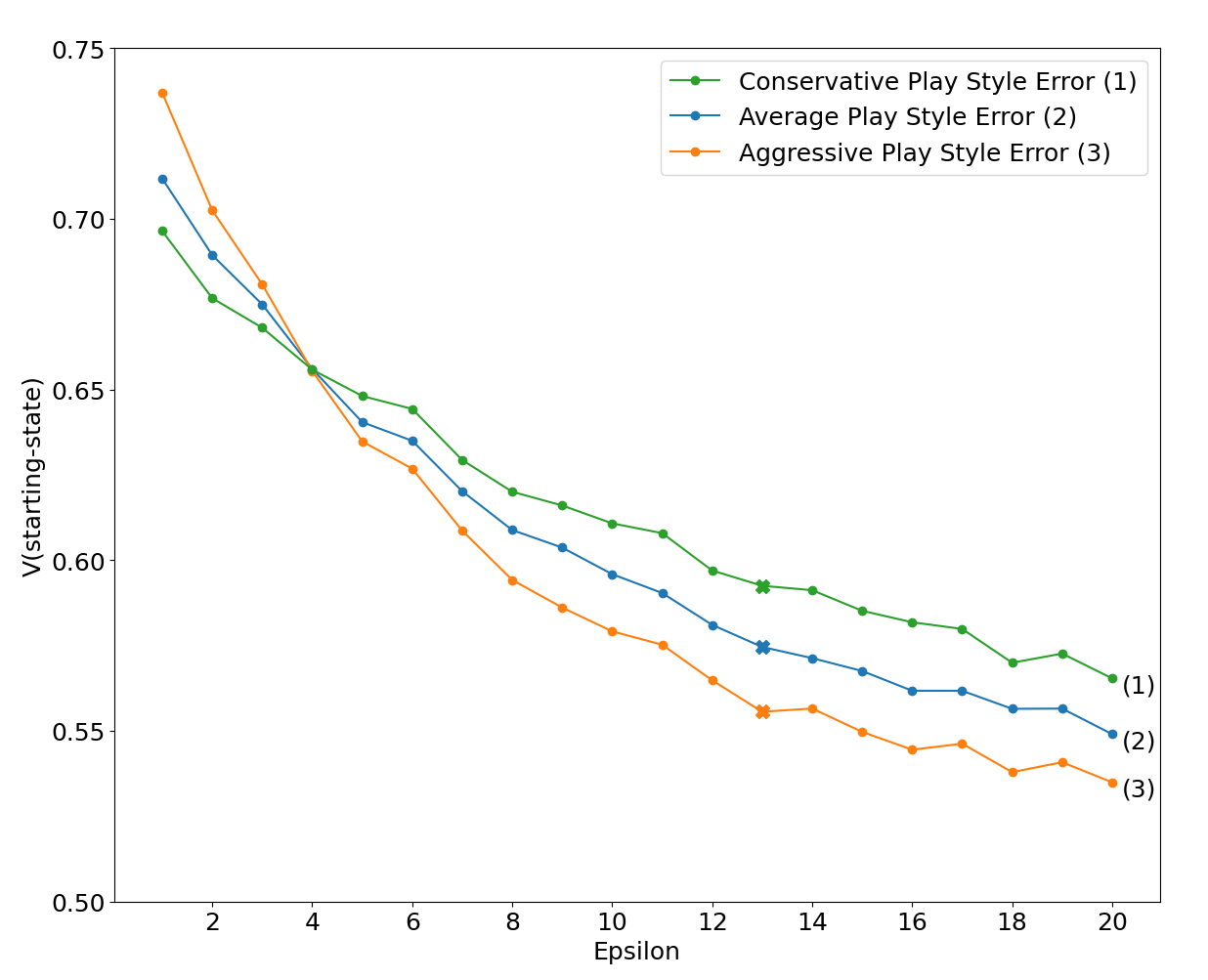}
    \caption{The relationship between value and error for players with average shot selection strategies compared to players who are aggressive or conservative.} % \tcyc{should we change "total error" to something else? } \cf{Value is measured as the weighted average for the value on both the serve and serve return states, denoted, $V(\text{starting-state)}$.}}
    \label{fig:mrp_playstyle}
\end{figure}

Figure \ref{fig:mrp_playstyle} depicts the relationship between value and execution error for an average player, conservative player and aggressive player. We can see that when execution error is small, it is advantageous to capitalize and be more aggressive. However, as execution error grows, continuing to be aggressive negatively impacts value relative to the average player. In contrast, the reverse relationship is true for conservative players. These results confirm general intuition that accurate players can afford to be more risk-seeking and take more aggressive shots. Interestingly, the switching point seems to occur around $\epsilon = 4$, which reflects a much more accurate shot maker than a player with average execution error ($\epsilon = 13$). Accordingly, our model suggests that a player with average execution error would be better off playing conservatively.

% \begin{figure}[H]
% \centering
% \begin{tabular}{cccc}
%     \includegraphics[width=0.30\textwidth, height=7.5cm]{Figures/playstyle_serve.png} &
%     \includegraphics[width=0.30\textwidth, height=7.5cm]{Figures/playstyle_servereturn.png}&
%     \includegraphics[width=0.30\textwidth, height=7.5cm]{Figures/playstyle_starting_state.png} \\
%     (a) Serve States & (b) Serve return states & (c) All states \\[6pt]
% \end{tabular}
% \caption{The relationship between value and error for players with average shot selection strategies compared to players that are aggressive or conservative. \tcyc{caption should clearly lay out what the three subfigs are. i could only infer from looking at y axes. also makes me wonder whether we need all 3, or having the "all states" one is enough and gives us a bigger fig to look at}}
% \label{fig:mrp_playstyle}
% \end{figure}

\subsection{Optimal value functions and policies in the presence of execution error}
\label{sec:MDP_Results}

\begin{figure}[H]
\centering
\begin{tabular}{cccc}
    \includegraphics[width=0.48\textwidth, height=4cm]{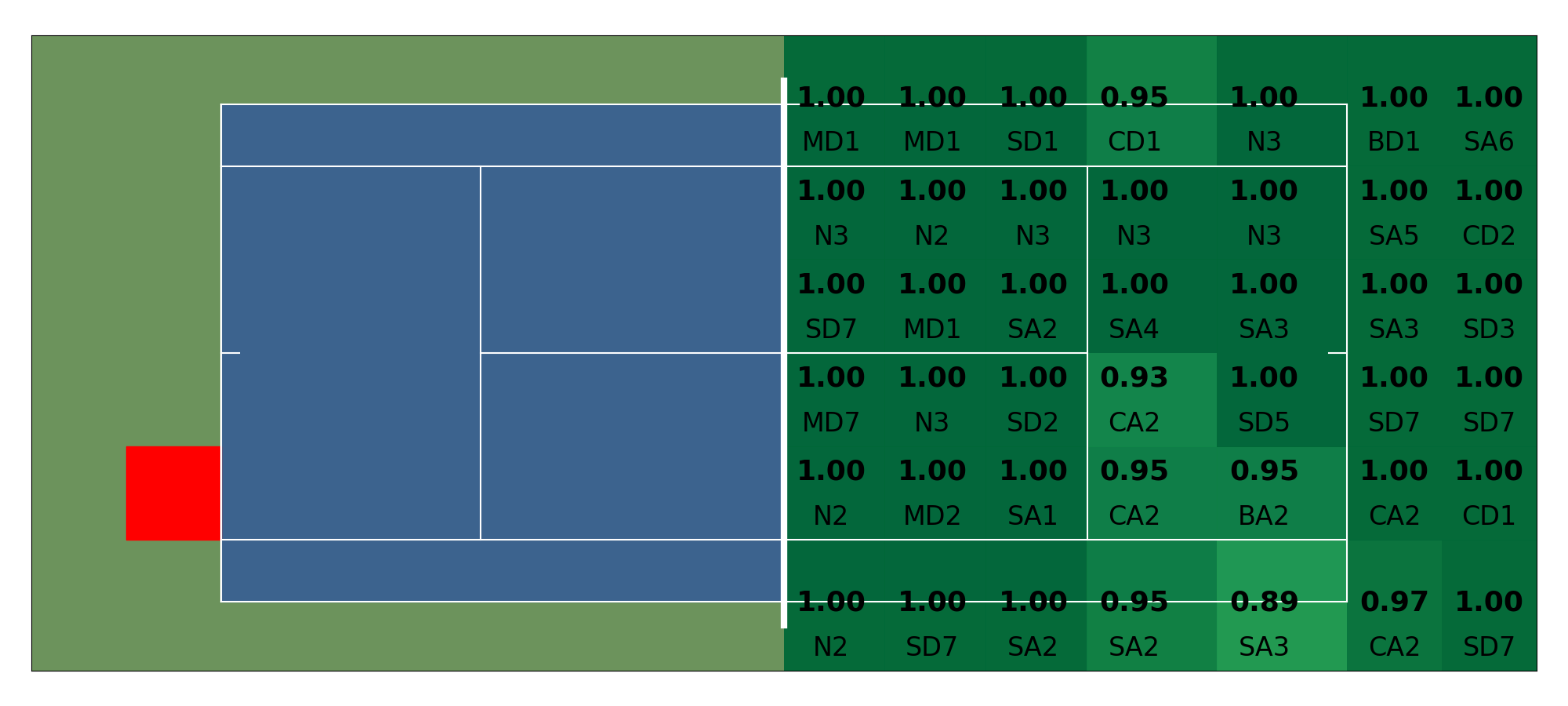} &
    \includegraphics[width=0.48\textwidth, height=4cm]{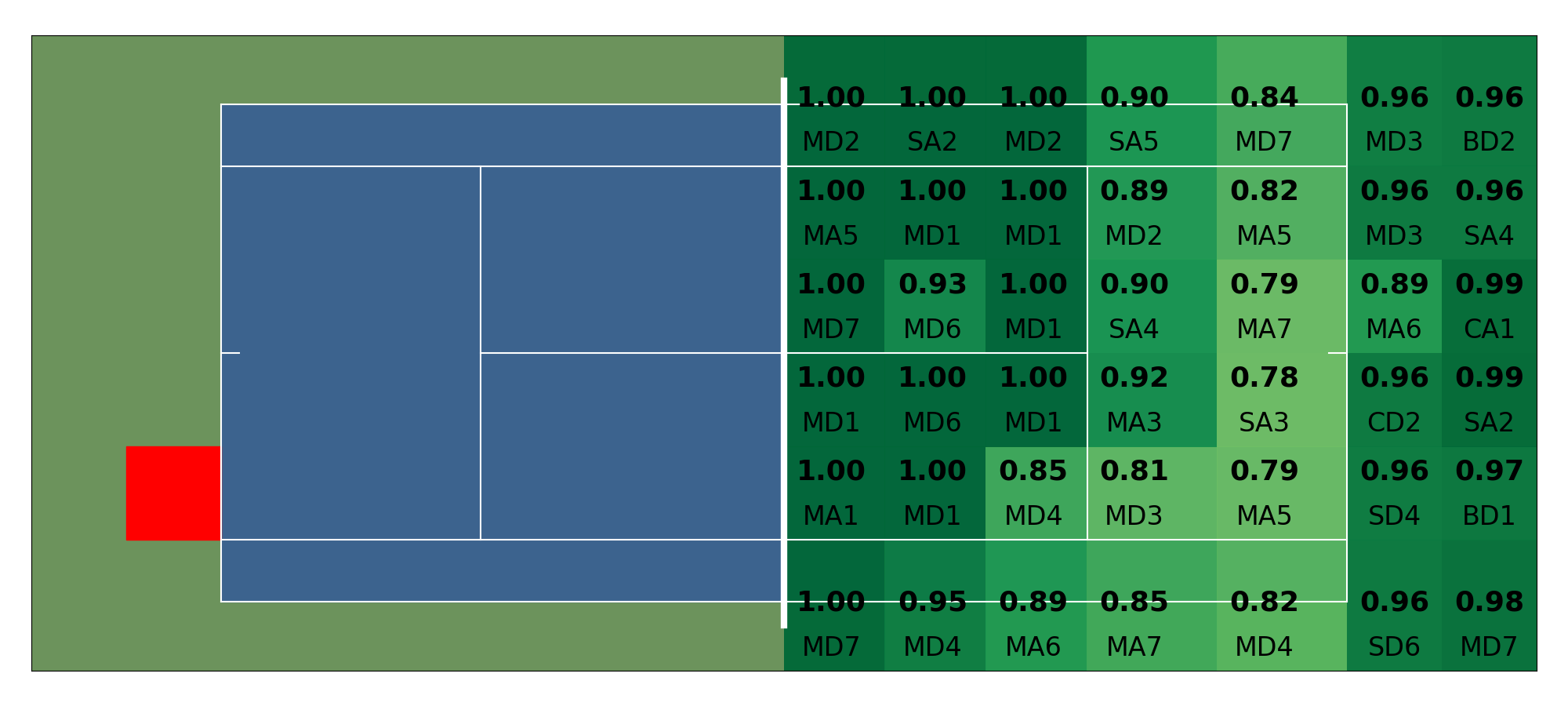} \\
    (a) Player A at deuce baseline with $\epsilon = 1$ & (b) Player A at deuce baseline with $\epsilon = 13$ \\[6pt]
\end{tabular}
\begin{tabular}{cccc}
    \includegraphics[width=0.48\textwidth, height=4cm]{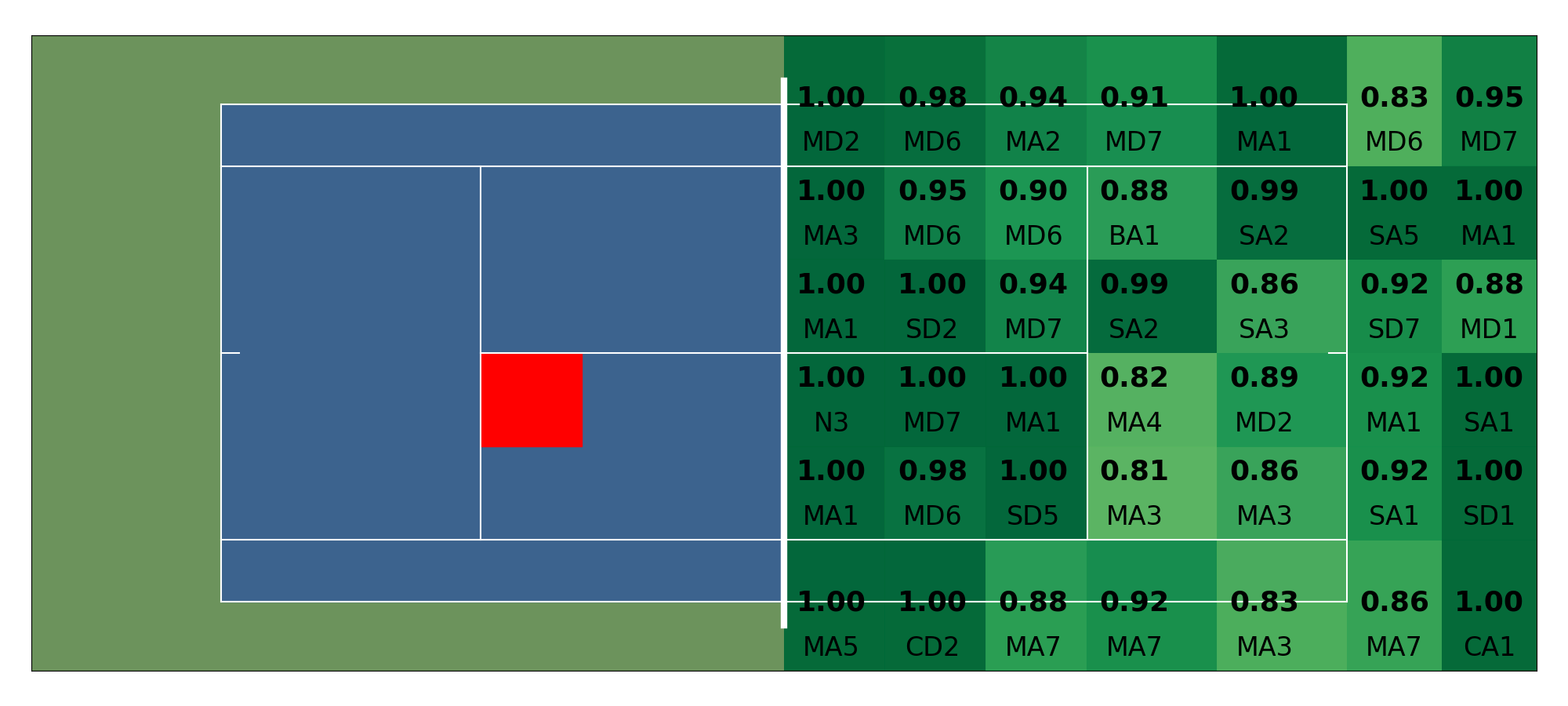} &
    \includegraphics[width=0.48\textwidth, height=4cm]{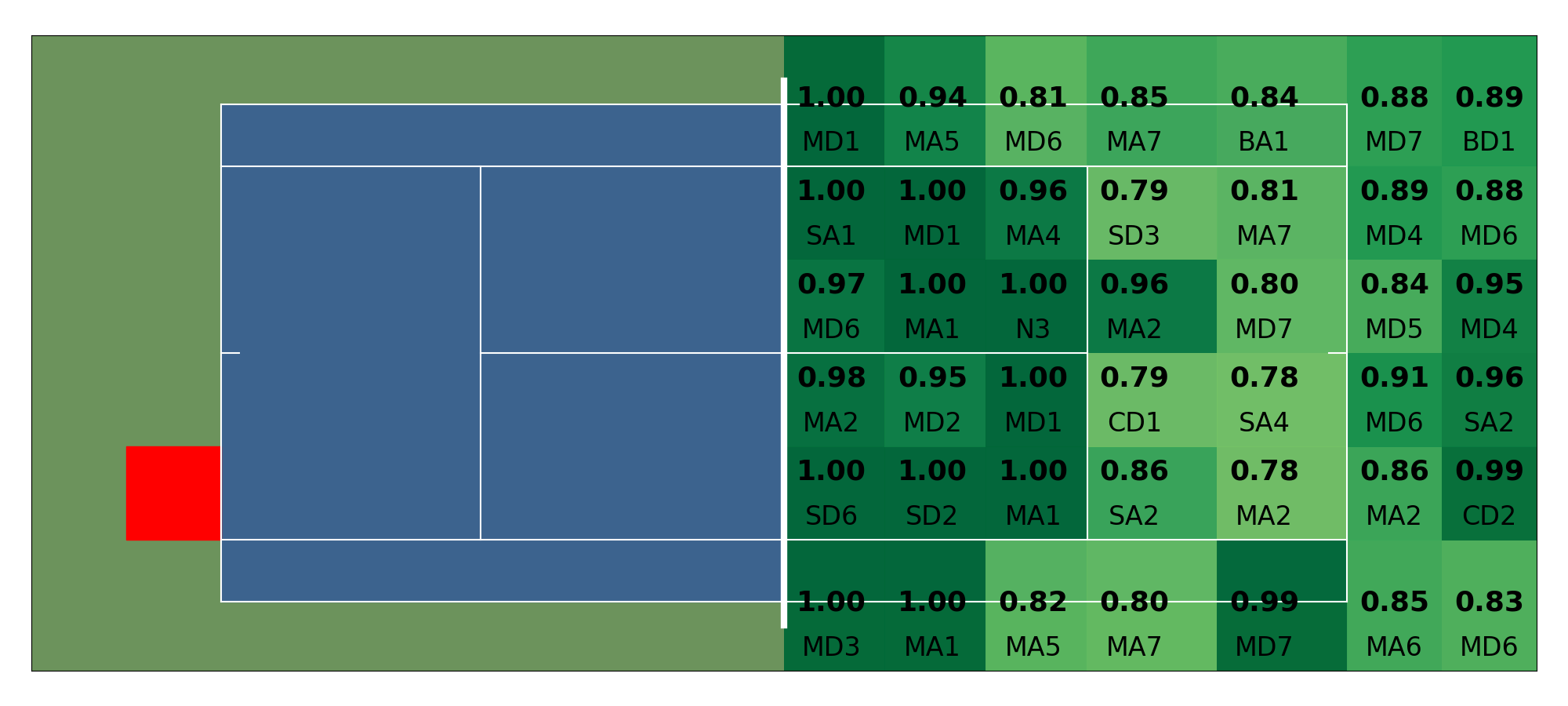} \\
    (c) Player A at deuce service line with $\epsilon = 1$ & (d) Player A at deuce baseline with $\epsilon = 20$\\[6pt]
\end{tabular}
\caption{The MDP value function for various rally states and various values of $\epsilon$. The number in each cell on the right side of the court indicates the value of that state, where the state is described as Player B standing in that cell, Player A standing in the red cell, and the shot being a rally shot being hit by Player A. The text in each cell indicates the optimal action for that state. %\tcyc{numbers in these figs are very hard to see. same with prev similar fig. need to address somehow.} \cf{will address this later}
}
\label{fig:mdp_Vs}
\end{figure}

Figure \ref{fig:mdp_Vs} shows the MDP results for several different rally shot states. As a reminder, $V^*(s)$ represents the probability that Player A wins the point starting from state $s$ when following the optimal shot selection strategy represented by $\pi^*$. The figure is read in the same manner as the MRP results in Figure \ref{fig:mrp_Vs}, with the number representing the win probability of the point, and the additional detail underneath indicating a corresponding optimal action. We see that there is an intuitive trend: value decreases as the execution error increases. However, the drop in value is much less acute than in the MRP case, for same change in $\epsilon$. For example, cross-court rally shots with $\epsilon = 20$ still have a value of around 89\%, in contrast to the 42\% for the same scenario in the MRP.

The reason for such high values for the MDP stems from the identification of optimal actions that lead to states with very high, nearly certain probabilities of winning the point. For example, the optimal policy will often choose shots that lead to the state where Player A is standing in a cell on the baseline and Player B is in cell 84 (deep on the deuce side behind the baseline). This is a rarely visited state in the MRP (and in real tennis) yet is a state that the MDP routinely visits because it is solving for an optimal shot selection strategy and certain actions lead to this state with a much higher probability. The act of optimizing leads to shot choices that are less realistic compared to what may be achievable in a real tennis match, while simultaneously exploiting small levels of sampling error in the transition probabilities from our simulation approach.

% \tcyc{i wouldn't say it is through the policy iteration alg. it is more that the transition dynamics lead us here. which suggests maybe there is a flaw in our transition matrix??} \cf{I dont think theres any flaws in the policy iteration or the transition matrix, just that the policy iteration will squeeze out any sampling error in the transition matrix to find ``shortest" paths that lead Player B out of position.} Symmetrically, unlike the MRP, the MDP never enters any state with a high probability of losing the point. Further discussion can be found in Appendix \ref{sec:AppendixC}

\begin{figure}[H]
\centering
\begin{tabular}{cccc}
    \includegraphics[width=12.5cm]{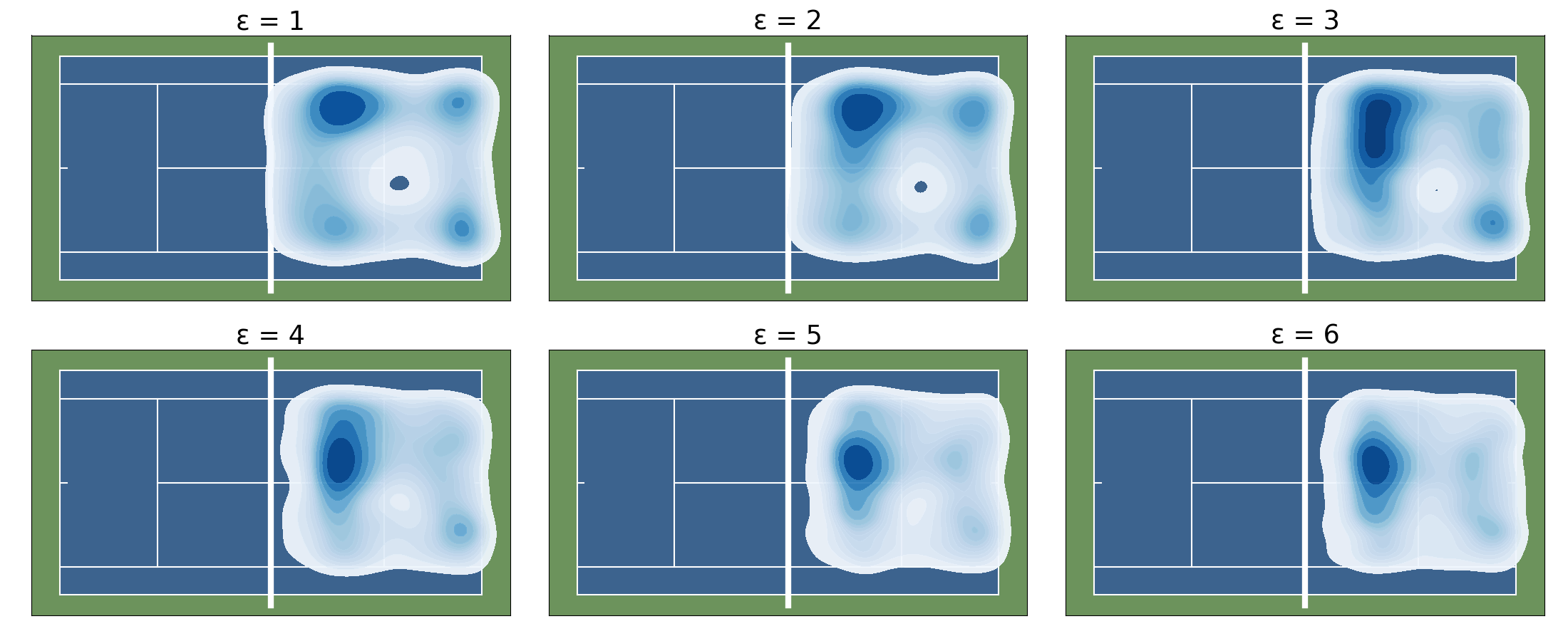} \\
    (a) Distribution of optimal actions over all states \\[6pt]
\end{tabular}
\begin{tabular}{cccc}
    \includegraphics[width=12.5cm]{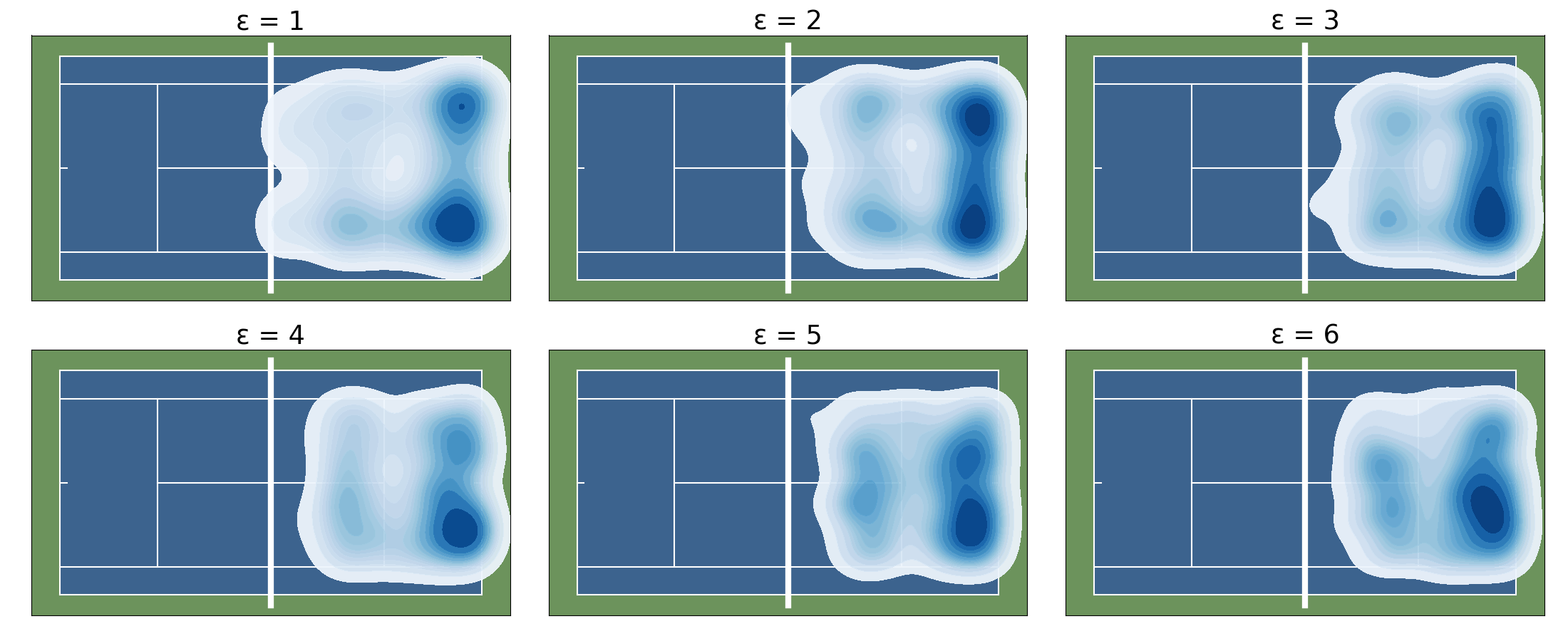} \\
    (b) Distribution of optimal actions for states with both players at the baseline \\[6pt]
\end{tabular}
\caption{The distribution of actions within the optimal policy for each level of error.}
\label{fig:fig:optimal_policy}
\end{figure}

Figure \ref{fig:fig:optimal_policy}(a) highlights how the distribution of optimal intentions over all states changes as a function of $\epsilon$.  The distribution is not symmetric due to the larger representation of simulated shots from right-handed players in the dataset. For $\epsilon = 1$, we see that the distribution peaks at intentions considered more aggressive (baseline corners, drop shots and short cross-court passes), with little to no probability mass on the conservative intentions. As $\epsilon$ increases from 1 to 6, the distribution shifts, moving probability mass away from the aggressive intentions at the baseline. For higher values of $\epsilon$, the distribution is similar to the $\epsilon = 6$ distribution, with a large majority of shots concentrated in the middle of the court. Figure \ref{fig:fig:optimal_policy}(b) depicts the same distribution but only over states where both players are standing behind the baseline. With low execution error, there is a high concentration of optimal actions aimed at the two baseline corners. This quickly changes as $\epsilon$ increases to favor shots aimed at the more central intentions, since aiming for the corners with larger execution error leads to more shots landing out of bounds.

Even though the MDP is behaving mathematically as it should, it still rests on the unrealistic premise that Player A can hit a ball at the optimal intention every shot. %\rs{[But doesn't the value function capture the long-run value?]} In a sense, the MDP is modeling the situation where Player A is going for a winner on every shot, without approach shots or other lower value shots that help set up the point. 
Moreover, our current state space does not capture the difficulty of the incoming shot that Player A needs to return, a factor that would undoubtedly impact the optimal policy and value. This, coupled with the optimization algorithm exploiting small levels of sampling error, likely results in an overly optimistic value function. We address this issue in the next section. 

%To address this issue, we developed the $n$-step lookahead approach to solving the MDP, instead of optimally solving the MDP. This results in more realistic and modest value functions, which we discuss in the next section. 

%This goes against regular tennis intuition that some shots, while having low value, are important and necessary to ``set up the winner." Rather, the MDP is modeling Player A ``going for the winner" on every shot, without ever having to use a set-up shot. This leads to higher values, however, less realistic ones. We discuss a tokenized approach to the MDP where Player A can only choose an optimal shot a select number of times within a point in Section \ref{sec:MDP-Regularized}.

\subsection{Varying the number of optimal actions in the policy}
\label{sec:MDP-nStep}

\begin{figure}[H]
    \centering
    \includegraphics[width=\textwidth]{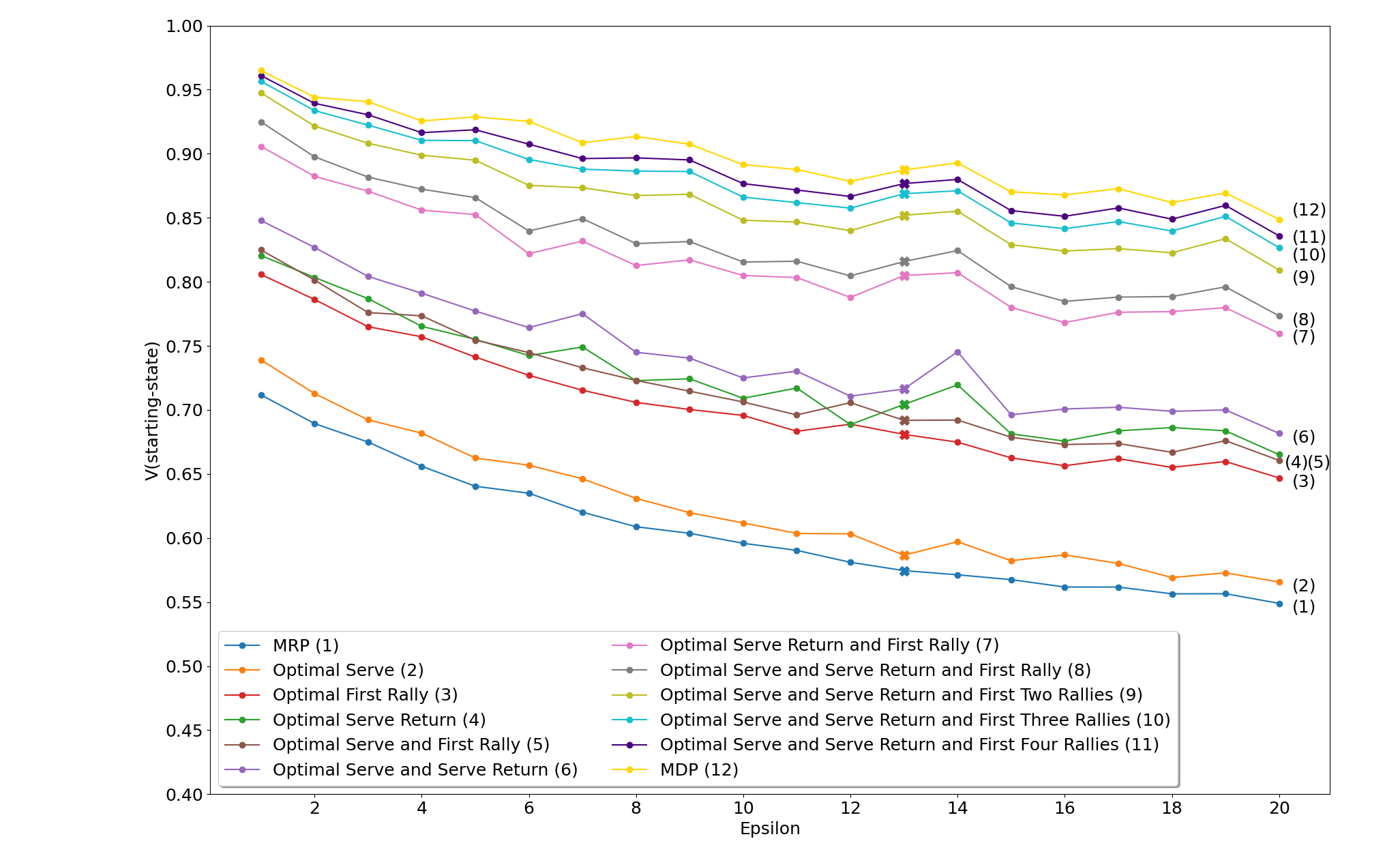}
    \caption{The relationship between value and error for different numbers and types of optimal shots.}% \cf{Value is measured as the weighted average for the value on both the serve and serve return states, denoted, $V(\text{starting-state)}$.}} 
    \label{fig:mdp_reg}
\end{figure}

In this final section, we use $\pi^n_\epsilon$ to quantify a more realistic value for optimal shot selection in the presence of execution error by allowing Player A to choose only a limited number of optimal actions. Figure \ref{fig:mdp_reg} highlights the relationship between value and execution error for the MRP with $\hat\pi$, the MDP with $\pi^*_\epsilon$, and a set of different $\pi^n_{\epsilon}$ policies modeling the following scenarios: %between one optimal action and five optimal actions at the start of the point. 

\begin{itemize}
    \item One optimal action (serve, serve return, or first rally shot)
    \item One optimal action (both serve and serve returns)
    \item Two optimal actions (serve and first rally shot or serve return and first rally shot) 
    \item Two optimal actions (both serve and serve returns, and first rally shot)
    \item Three optimal actions (both serve and serve returns, and first two rally shots)
    \item Four optimal actions (both serve and serve returns, and first three rally shots)
    \item Five optimal actions (both serve and serve returns, and first four rally shots)
\end{itemize}

In accordance with theory, we see that the more optimal actions available to Player A, the higher the value. The figure also illustrates diminishing returns in value as the number of optimal actions increases. Furthermore, the majority of points in tennis last fewer than eight shots (four per player), so the advantage of having five optimally selected shots is quite close to the value of the MDP, which represents optimal shot selection on every shot in every point. This result is also consistent with our observation that policy iteration generally terminated after 6-8 iterations.

If Player A only has one shot to optimize, the biggest gain in value is achieved by optimizing the serve return. This is due to the server advantage in tennis and holds true for all values of execution error. This large gain in value from optimizing serve returns can be interpreted as the particular shot players should put the most effort into hitting well, as it has the highest return on investment. In fact, optimizing the service return actually has comparable value to optimizing both serves and first rally shots for small epsilon, and actually slightly higher value as epsilon increases. 

%This finding motivates the importance of practicing the serve return in training schedules.  

Going deeper, if a player with average execution error ($\epsilon = 13$) can optimize their serve return, the probability of winning the point goes from 0.575 to 0.704. This increase in value is roughly equivalent to employing the policy $\hat\pi$ (i.e., no optimization of the intentions) but with an execution error of $\epsilon = 1$. In other words, being able to choose a single optimal intention for serve returns and then choosing actions according to the intention distribution, with average execution error for all shots, is equivalent to playing according to the intention distribution with perfect execution for all shots. Although this increase in value from a single optimized intention may seem high, recall that our models are based on shots generated by VON CRAMM, which is trained on data from professional matches. Being able to hit an optimal return shot against a professional player is not an easy task. Similarly, average execution error should be interpreted as the average of top professional players in the world. Nevertheless, these results provide insight into the untangling of the contribution of intention versus execution to value generation. %\tcyc{[how is this? i moved and expanded it here. seemed a bit out of place above]} \cf{[I like it]}

%\rs{[Does this mean pros do not employ optimal serve-return?  Probably such data is not available but would be cool to do pro-specific analysis.  For example, compare a pro's shot-selection to the optimal policy.  I am guessing top players would be near-optimal.  And perhaps one can understand the difference between the top players and others.  Quite possible that every decent player employs near-optimal shot-selection but the top players have a lower epsilon.  Also possible that some players have a higher epsilon than others but are better because they adapt to their epsilon (i.e., find the right balance between being agressive/conservative).]} \cf{[tried to add to next steps]} \tcyc{[we may want to go back and reiterate that this von cramm simulated data is modeled on top pros, right? and therefore average epsilon is average error of top pros, so already quite good. there may be a misconception that when we say "average" people intuitively think of more recreational average]} \cf{[Thats fair- I don't think we should say "top" pros, but regular pro, since von cramm is built on all rounds of the Australian open. Added a sentence to sec 5.1]}

Overall, these results help untangle the contribution of intention and execution error to value. The probability of winning the point can be increased via both choosing where to hit the ball and executing the shot well. Given a player with average execution error, optimizing one action, either the serve return or first rally shot, has a much bigger impact on the win probability compared to reducing execution error. One would need perfect execution to generate similar value, which seems more difficult than selecting the right shot in a given state.

%Moving horizontally from the optimal serve line until we intersect the MRP line, we see that if a player chooses to not optimize the serve, they need to have an error level of $\epsilon = 11$ to maintain roughly the same value. Similarly, if the player chose to optimize only the first rally shot, their probability of winning the point jumps to 0.681. In order to maintain this value with no optimal epochs, the player would need an error level of $\epsilon = 2.5$, something that is highly unlikely.

% By construction, the one-step lookahead policy models the situation where Player A can choose an optimal action in the current decision epoch and then follows the empirical policy $\pi^e$ for the remaining epochs. \tcyc{some of this text should probably go into the methods section above. above i guess we should define the $n$ step lookahead since the results have $n$ more than 1. i suggest we define 1 first to make things concrete and then generalize to $n$. }

% For example, if we want to model the behavior of only being able to choose one optimal action per point, we can choose an optimal serve and then continue with the empirical policy. Alternatively, if we want to measure the value of two optimal epochs, we can try being optimal at our first and second rally shot, and non-optimal for all remaining shots. \tcyc{have to be careful since 2 opt shots is not the same as the first two shots being opt. we could discuss this wrt expanding the state space and how we wish not to do that.}

% \begin{figure}[H]
% \centering
% \begin{tabular}{cccc}
%     \includegraphics[width=0.31\textwidth, height=7.5cm]{Figures/mdp_reg_serve.png} &
%     \includegraphics[width=0.31\textwidth, height=7.5cm]{Figures/mdp_reg_serve_return.png} &
%     \includegraphics[width=0.31\textwidth, height=7.5cm]{Figures/mdp_reg_starting_state.png} \\
%     (a)  & (b) & (c) \\[6pt]
% \end{tabular}
% \caption{The relationship between value and error with various number and types of optimal epochs. \tcyc{the fig might be easier to interpret if we got rid of the legend and had the labels on the RHS of plot with arrows pointing to the left indicating which line? edit: i added the itemize list, which may help. but still might need to make the graph easier to read, esp in black and white}}
% \label{fig:mdp_reg}
% \end{figure}

%\tcyc{STOPPED HERE}

\section{Conclusion}
\label{sec:Conclusion}

\noindent
In this paper, we develop a modeling framework comprising Markov reward and Markov decision processes to investigate the impact of execution error on the likelihood of winning and optimal strategy in sport. We apply this framework to tennis and derive several tennis-specific insights. Our analysis shows that erring on backhand shots is more costly than on forehand shots and that execution error on serves has a more limited impact on the probability of winning the point. We demonstrate that larger execution error leads to more conservative shot selection strategies, and identify the break-even point in execution error between aggressive and conservative strategies being optimal. Lastly, we show that optimal shot selection on serve returns is more valuable than on any other shot, over all values of execution error. Our results give insight into the types of shots in which player should exert more effort, and for a given shot, the change in the probability of winning the point when choosing an optimal shot versus reducing execution error. %\tcyc{i think it may be tricky to get across the diff between optimal a and reducing eps in abstract/conclusion, when there isn't much space yet to explain everything. also wondering if we should replace "value" with "prob of winning point" in abstract/conclusion to be more concrete and less conceptual, for the same reason} \cf{agree - makes sense. Tried to change that here in conclusion.}

%by exploring the impact of certain optimal epochs, we motivate the importance of practicing and exerting the most energy on serve returns, as this has the highest incremental gain in value. 

Our modeling framework can be expanded within tennis in several ways. First, the state space could be enriched to include information such as the incoming shot's difficulty and the current score. Second, whereas we currently define execution error as a deviation from the desired landing location of the ball, this could be expanded to include deviations from the desired ball speed and spin. Third, if one has access to ball and player tracking data for specific players, individualized value functions and optimal policies could be determined. This data can also be used to determine the precise $\epsilon$ and policies specific players employ.

Finally, we note that our modeling framework is general and broadly applicable beyond tennis. Most readily, other racket and net sports such as badminton and table tennis would be straightforward to analyze. Other sports that have very clear intended targets would also be suitable, such as curling, bowling and darts. Lastly, we highlight the opportunity of expanding the MDP model to allow Player B to also choose optimal shots, instead of just following the intention distribution. This approach may be of special interest to game theorists.

\section*{Acknowledgments} 
We would like to thank Raghav Singal and Martin Puterman for helpful suggestions and comments that improved our paper. This work was partially funded by the Connaught Global Challenge Award (2021-2021).

\bibliographystyle{elsarticle-harv} 
\bibliography{covidbib}

\appendix
\section{Determining the error level of an average player}
\label{sec:AppendixA}

Table \ref{table:averageError} shows the probabilities of winning,\textbf{ erroring }or keeping shots in player for different values of $\epsilon$. In order to determine which $\epsilon$ value most closely follows the average player, we look to see where these values begin to converge with the empirical results. This is seen most closely with $\epsilon=13$.

\begin{table}[H]
\centering
	\captionsetup{singlelinecheck = false, justification=raggedright}
	\captions{The probability of winning, erroring or keeping shots in-play for each level of $\epsilon$, compared to the empirical probabilities for the average player. The closest relationship exists with $\epsilon=13$.}
	\label{table:averageError}
\renewcommand{\arraystretch}{1}% Tighter
\begin{tabular}{|c|c|c|c|}
\hline
\rowcolor[HTML]{C0C0C0} 
\cellcolor[HTML]{9B9B9B}                 & \textbf{P(win)} & \textbf{P(error)} & \textbf{P(in-play)} \\ \hline
\cellcolor[HTML]{C0C0C0}$\epsilon=1$     & 14.4            & 0.4               & 85.2                \\ \hline
\rowcolor[HTML]{EFEFEF} 
\cellcolor[HTML]{C0C0C0}$\epsilon=2$     & 14.4            & 1.8               & 83.8                \\ \hline
\cellcolor[HTML]{C0C0C0}$\epsilon=3$     & 14.0            & 3.2               & 82.8                \\ \hline
\rowcolor[HTML]{EFEFEF} 
\cellcolor[HTML]{C0C0C0}$\epsilon=4$     & 14.1            & 4.8               & 81.1                \\ \hline
\cellcolor[HTML]{C0C0C0}$\epsilon=5$     & 13.9            & 6.1               & 80.0                \\ \hline
\rowcolor[HTML]{EFEFEF} 
\cellcolor[HTML]{C0C0C0}$\epsilon=6$     & 13.6            & 7.3               & 79.1                \\ \hline
\cellcolor[HTML]{C0C0C0}$\epsilon=7$     & 13.6            & 8.2               & 78.2                \\ \hline
\rowcolor[HTML]{EFEFEF} 
\cellcolor[HTML]{C0C0C0}$\epsilon=8$     & 13.4            & 9.0               & 77.6                \\ \hline
\cellcolor[HTML]{C0C0C0}$\epsilon=9$     & 13.6            & 10.1              & 76.3                \\ \hline
\rowcolor[HTML]{EFEFEF} 
\cellcolor[HTML]{C0C0C0}$\epsilon=10$     & 13.2            & 10.8              & 76.0                \\ \hline
\cellcolor[HTML]{C0C0C0}$\epsilon=11$     & 13.3            & 11.6              & 75.1                \\ \hline
\rowcolor[HTML]{EFEFEF} 
\cellcolor[HTML]{C0C0C0}$\epsilon=12$     & 13.2            & 12.3              & 74.5                \\ \hline
\rowcolor[HTML]{FCFF2F} 
\cellcolor[HTML]{C0C0C0}$\epsilon=13$     & 13.0            & 13.0              & 74.0                \\ \hline
\rowcolor[HTML]{EFEFEF} 
\cellcolor[HTML]{C0C0C0}$\epsilon=14$     & 13.0            & 13.6              & 73.4                \\ \hline
\cellcolor[HTML]{C0C0C0}$\epsilon=15$     & 12.9            & 14.2              & 72.9                \\ \hline
\rowcolor[HTML]{EFEFEF} 
\cellcolor[HTML]{C0C0C0}$\epsilon=16$     & 13.0            & 15.0              & 72.0                \\ \hline
\cellcolor[HTML]{C0C0C0}$\epsilon=17$     & 12.9            & 15.4              & 71.7                \\ \hline
\rowcolor[HTML]{EFEFEF} 
\cellcolor[HTML]{C0C0C0}$\epsilon=18$     & 13.0            & 16.1              & 70.9                \\ \hline
\cellcolor[HTML]{C0C0C0}$\epsilon=19$     & 12.8            & 16.8              & 70.4                \\ \hline
\rowcolor[HTML]{EFEFEF} 
\cellcolor[HTML]{C0C0C0}$\epsilon=20$     & 12.5            & 17.1              & 70.4                \\ \hline
\rowcolor[HTML]{9B9B9B} 
\multicolumn{4}{|c|}{\cellcolor[HTML]{9B9B9B}}                                                       \\ \hline
\rowcolor[HTML]{FCFF2F} 
\cellcolor[HTML]{C0C0C0}\textbf{Average Player} & 13.0            & 13.4              & 73.6                \\ \hline
\end{tabular}
\end{table}

\section{Explaining the MRP value function for ad vs. deuce rally shots}
\label{sec:AppendixB}

\begin{figure}[H]
\centering
\begin{tabular}{cccc}
    \includegraphics[width=0.45\textwidth]{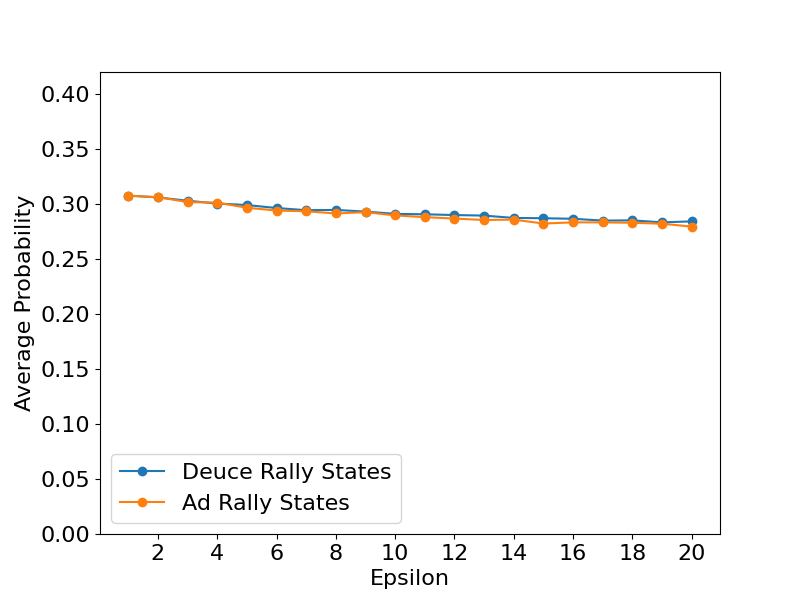} &
    \includegraphics[width=0.45\textwidth]{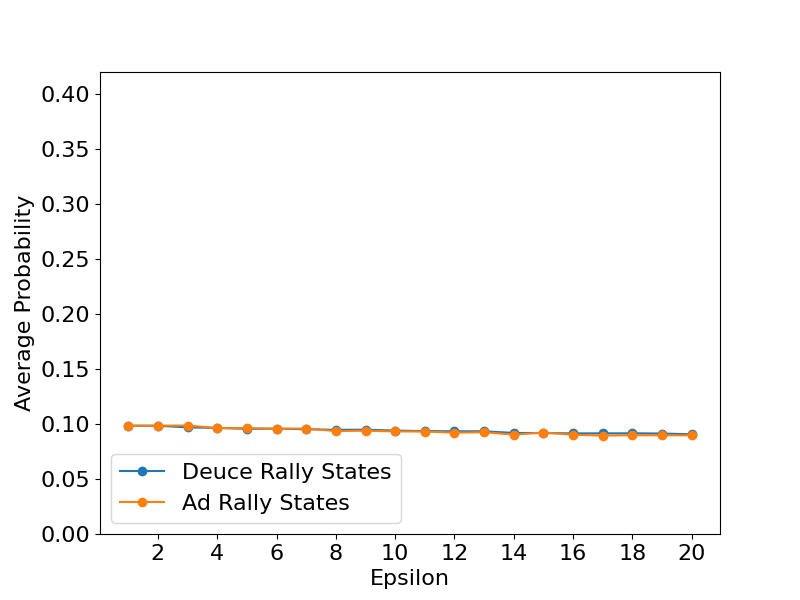} \\
    (a) Transitioning into Player A's shot winning & (b) Transitioning into Player B's shot erring \\[6pt]
\end{tabular}
\begin{tabular}{cccc}
    \includegraphics[width=0.45\textwidth]{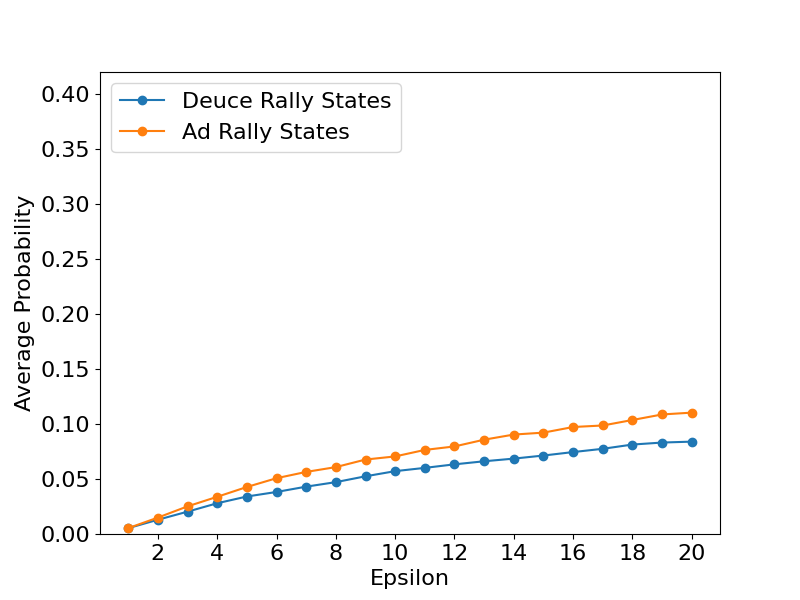} &
    \includegraphics[width=0.45\textwidth]{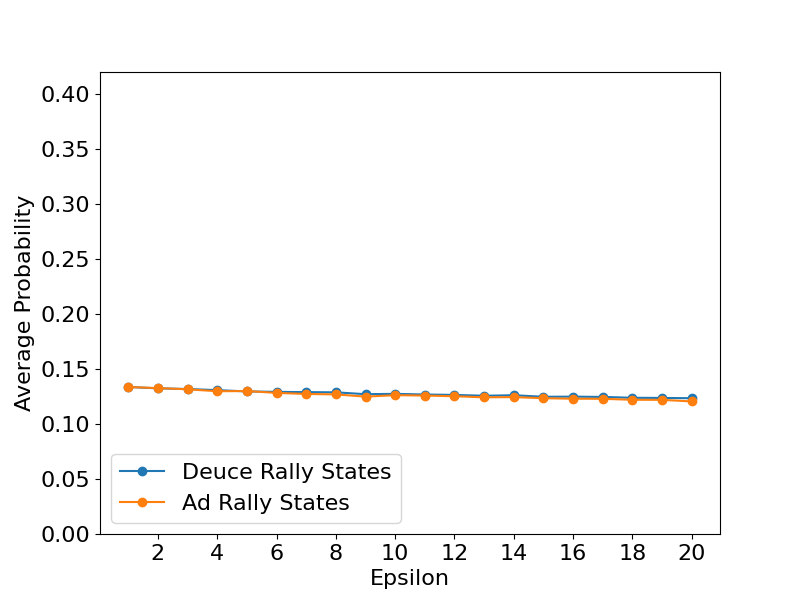} \\
    (c) Transitioning into Player A's shot erring & (d) Transitioning into Player B's shot winng \\[6pt]
\end{tabular}
\caption{The average probabilities for Player A and B winning and erring the shot from the ad and deuce rally states at various levels of $\epsilon$.}
\label{fig:ad_deuce}
\end{figure}

Figure \ref{fig:ad_deuce} shows the average transition probabilities for Player A and B winning and erring the shot from the ad and deuce rally states at various levels of $\epsilon$. We note that the trends between ad and rally states are indistinguishable except in Figure \ref{fig:ad_deuce}(c), the probability of Player A's shot erring. As $\epsilon$ increases, this probability grows at faster rate for ad rallies compared to deuce rallies. Due to the high prevalence of right-handed players and noting that backhanded shots are generally harder to hit, this result is intuitive.

% \section{Explaining the high MDP value function}
% \label{sec:AppendixC}

% \cf{Need to add the figures here - will do later.}

\end{document}